\documentclass[11pt]{amsart}
\usepackage{amsmath,amssymb,amsfonts,url,mathptmx,txfonts}
\usepackage{color}
\usepackage{appendix}
\allowdisplaybreaks[4]
\numberwithin{equation}{section}

\newtheorem{thm}{Theorem}[section]

\usepackage[marginal]{footmisc}

\usepackage{amsthm,a0size,bm,indentfirst,wasysym}

\def\XXint#1#2#3{{\setbox0=\hbox{$#1{#2#3}{\int}$ }
		\vcenter{\hbox{$#2#3$ }}\kern-.6\wd0}}

\newtheorem{lem}[thm]{Lemma}

\newtheorem{pro}[thm]{Proposition}
\newtheorem{rem}[thm]{Remark}

\usepackage{amsmath}

\allowdisplaybreaks[4]

\renewcommand{\a }{\alpha }
\renewcommand{\b }{\beta }
\renewcommand{\d}{\delta }
\newcommand{\D }{\Delta }

\newcommand{\e }{\varepsilon }
\newcommand{\g }{\gamma}
\renewcommand{\i }{\iota}
\newcommand{\G }{\Gamma }
\renewcommand{\l }{\lambda }

\newcommand{\n }{\nabla }

\newcommand{\s }{\sigma }

\renewcommand{\o }{\omega }

\newcommand{\ov}{\overline}
\newcommand{\intbar}{\mathop{\int\makebox(-13.5,0){\rule[4pt]{.7em}{0.3pt}}%
\kern-6pt}\nolimits}

\newcommand{\be}{\begin{equation}}
\newcommand{\ee}{\end{equation}}
\newcommand{\bes}{\begin{equation*}}
\newcommand{\ees}{\end{equation*}}
\newcommand{\ba}{\begin{eqnarray}}
\newcommand{\ea}{\end{eqnarray}}
\newcommand{\bas}{\begin{eqnarray*}}
\newcommand{\eas}{\end{eqnarray*}}
\newenvironment{pf}{\noindent{\sc Proof}.\enspace}{\rule{2mm}{2mm}\medskip}
\newenvironment{pfn}{\noindent{\sc \bf Proof}}{\rule{2mm}{2mm}\medskip}

\newtheorem{remark}{Remark}[section]
\newcommand{\R}{\mathbb{R}}

\newcommand{\N}{\mathbb{N}}
\renewcommand{\o }{\omega }

\begin{document}

\title[Nirenberg problem ]{The Nirenberg problem on   half spheres: \\
A bubbling off analysis }

\author{Mohameden Ahmedou}

\author{Mohamed Ben Ayed }

\address{ Mohameden Ahmedou\\
Department of Mathematics \\
	Giessen University \\
	Arndtstrasse 2, 35392, Giessen, Germany }
\email{Mohameden.Ahmedou@math.uni-giessen.de}

\address{  Mohamed Ben Ayed, Department of Mathematics, College of Science, Qassim University, Buraidah 51452, Saudi Arabia \\ \& 
Universit\'e de Sfax, Facult\'e des Sciences de sfax, Route de Soukra, Sfax, BP. 1171, 3000 Tunisia } 
\email{M.BenAyed@qu.edu.sa \& Mohamed.Benayed@fss.rnu.tn  }



 \maketitle


\centerline{   \emph{dedicated to the  memory of Prof. Antonio Ambrosetti }}

\bigskip

\begin{center}
{\bf Abstract}
\end{center}
In this paper we perform a refined blow up analysis of finite energy approximated solutions to a Nirenberg type  problem on half spheres. The later   consists of prescribing, under minimal boundary conditions,  the scalar curvature to be a given function. In particular we give a precise location of blow up points and blow up rates. Such an analysis shows that the blow up picture of the Nirenberg problem on half spheres is far more complicated that in the case of closed spheres. Indeed besides the combination of interior and boundary blow ups, there are non simple blow up points for subcritical solutions having zero or  nonzero weak limit. The formation  of such  non simple blow ups is governed by a vortex problem, unveiling an unexpected connection with  Euler equations in fluid dynamic and mean fields type equations in mathematical physics.

\begin{center}


\centerline{{\bf AMS subject classification:}  58J05, 35A01, 58E05.}

\end{center}

\tableofcontents

\section{Introduction and statement of main results}

The Nirenberg problem  in conformal geometry  has a long history. Indeed this  problem  goes back to the following question raised by Louis Nirenberg in the academic year  1969-1970: Given a smooth function $K$  on the standard  sphere $\mathbb{S}^n$ endowed with its standard metric $g_0$,  does there exist a metric $g$,  in the  conformal class of  $g_0$,  whose scalar curvature is given by the function $K$?
On the two dimensional sphere $\mathbb{S}^2$ this amounts to solve an elliptic PDE involving exponential nonlinearity while on spheres $\mathbb{S}^n$ of dimension $n \geq 3$, using the transformation law of the scalar curvature under conformal change of metric, one sees that the Nirenberg problem is  equivalent to solving the following nonlinear partial differential equation involving the critical Sobolev exponent:
$$
(\mathcal{NP}) \qquad L_{g_0} u \, = \, K u^{(n+2)/(n-2)}, \qquad u > 0,
$$
where  $ L_{g_{0}}:= -\D_{g_{0}}  \, + \, {n(n-2)}/{4}  $ denotes the conformal Laplacian.\\
The Nirenberg problem  has attracted a lot of attention in the last half century. See \cite{Au76, AH91,  BC2,  Bahri-Invariant, BCCH, BN83, CY, CGY, CL, CL1, CL2, Chen-Xu, Hebey, KW1, yyli1, yyli2, SZ} and the references therein.\\
Due to Kazdan-Warner obstructions \cite{KW1, BEZ1} the above PDE is not always solvable and the corresponding Euler-Lagrange functional lacks compactness.
 One way to overcome the lack of compactness, which goes back  to  Yamabe \cite{Yamabe} and R. Schoen \cite{Schoen} (see also \cite{DHR}) consists of considering the following subcritical approximation
$$
(\mathcal{NP}_{\e}) \qquad  L_{g_0} u \, = \, K u^{((n+2)/(n-2))- \e}, \qquad u > 0
$$
whose Euler-Lagrange functional satisfies the Palais-Smale condition. One then studies the blow up behavior when $\e$ goes to zero. It follows from such a refined  analysis that blow ups occur at critical point with nonpositive
Laplacian. Moreover  it turns out that the  blow up scenario depends strongly on the dimension and on the behavior of the function $K$ around its critical points. Indeed  under the condition that $\D K \neq 0 $  at the  critical points of $K$, one has  that,  in dimension $3$, solutions of the above approximation $ (\mathcal{N P}_{\e})$ could  only develop one single bubble,  see please  \cite{BC2, CGY}.  In  dimension $4$ there could be multiple blow ups  but tuples of blow up points   have
to  satisfy a balancing condition, see please \cite{BCCH, yyli2},  while in dimensions $n \geq 5$ every tuple of  distinct critical points of $K$  having negative Laplacian can be realized as  the blow up set of a blowing up solution of the approximated PDE, see please \cite{Bahri-Invariant, MM, MM19}. Furthermore all  blow up poins of finite energy blowing up solutions are  \emph{isolated simple} in the sense that  around every blow up point there is a ball which does not contain any other blow up point
and the Dirichlet-Energy of the blow up solutions $u_{\e}$ in a shrinking neighborhood around the blow up  point tends as $\e \to 0$ to the energy of \emph{one bubble} concentrating at this point. See please \cite{MM, MM19}.
\medskip

In this paper we consider  a version of the Nirenberg problem on standard half spheres $(\mathbb{S}^n_+,g)$. Namely we prescribe simultaneously   the scalar curvature to be a positive function $ K \in C^3(\mathbb{S}^n_+)$ and the boundary mean curvature to be zero. This amounts to solve the following boundary value problem
\begin{equation}
(\mathcal{P}) \quad
\begin{cases}
   -\D_{g_{0}} u \, + \, \frac{n(n-2)}{4} u \,  =  K \, u^{{(n+2)}/{(n-2)}},\,  u > 0 & \mbox{ in } \mathbb{S}^n_+, \\
  \frac{\partial u}{\partial \nu }\, =\, 0  & \mbox{on } \partial \mathbb{S}^n_+.
\end{cases}
\end{equation}

 This problem has been studied  on half spheres of dimensions $ n = 3, 4 $. See  the papers \cite{BEOA, BEO, BGO, BOA, DMOA, yyli, LL} and the references therein.
Like  the case of the Nirenberg problem on spheres, there are obstructions to the existence of solutions to $(\mathcal{P})$ and the corresponding variational problem is not compact. In order  to recover compactness  one considers    the  following subcritical approximation
 \begin{equation}
(\mathcal{P}_{\e}) \quad
\begin{cases}
   -\D_{g_{0}} u \, + \, \frac{n(n-2)}{4} u \,  =  K \, u^{((n+2)/(n-2)) - \e}, \,  u > 0 & \mbox{ in } \mathbb{S}^n_+, \\
  \frac{\partial u}{\partial \nu }\, =\, 0  & \mbox{on } \partial \mathbb{S}^n_+.
\end{cases}
\end{equation}
Regarding  the behavior of a sequence of  energy bounded  solutions $u_{\e}$ of $(\mathcal{P}_{\e})$, it follows from the concentration compactness principle that  either the $||u_{\e}||_{L^{\infty}}$ remains uniformly bounded  or  $u^{{2n}/{(n-2)}}_{\e}  \mathcal{L}^n$ (where $ \mathcal{L}^n $ denotes the Lebesgue measure) converges in the sense of measures  to a sum of Dirac masses, some of them are located  in the interior and the others
are  on the boundary. Moreover  it follows from the blow up analysis of $(\mathcal{P}_{\e}) $ that the interior points are  critical points of $K$  with nonpositive
  Laplacian and the boundary points are critical points of $K_1$ the restriction of $K$ on the boundary and satisfying that ${\partial_{\nu} K} \geq  0$. See \cite{BEO, BOA, DMOA}.\\
Furthermore, under the  non degeneracy assumption that $\D K \neq 0$ at interior critical points of $K$ and that ${\partial_{\nu} K} \neq 0$ at critical points of $K_1$, we have that    in the dimension $n=3$  multiple bubbling may occur but all  blow up points are \emph{isolated simple}, see \cite{DMOA, yyli}. Moreover   under additional condition on $K_1$ it has been proved in \cite{BGO} that in dimension $4$ all blow up  points are isolated simple. More surprisingly and in contrast with  the case of closed spheres, the Nirenberg  problem on half spheres may have \emph{non simple blow up points} for finite energy bubbling solutions of  $(\mathcal{P}_\e)$ see \cite{AB20b}.
\medskip

In this paper we perform a systematic asymptotic analysis, as $\e$ goes to zero,  of  finite energy blowing up solutions
 of $ (\mathcal{P}_{\e})$. Such an analysis is performed under the following non degeneracy conditions: \\
$\textbf{(H1)}$ The critical  points $y$'s of $K$ in $\mathbb{S}^n_+$ are non degenerate and  satisfy $\D K(y) \neq 0$.\\
$\textbf{(H2)}$ The critical points $z$'s of  the restriction of $K$ on the boundary $K_1:= K_{| \partial \mathbb{S}^n_+}$ are non degenerate and  satisfy $ \partial_{\nu} K (z)\neq 0$.\\
We first consider the case where the sequence of energy bounded solutions of $(\mathcal{P}_\e)$ has a zero weak limit  $u_\e \rightharpoonup 0$. In this situation the description of the blow up picture is as follows:
\begin{thm}\label{th:t1}
 Let $n\geq 5$ and  $0 < K \in C^3(\ov{\mathbb{S}^n_+})$ be a positive function  satisfying  the assumptions  $(H1), \, (H2)$ and let $(u_\e)$ be a sequence of energy bounded solutions of $(\mathcal{P}_\e)$ with $u_\e \rightharpoonup 0$. Then $u_\e$  blows up and decomposes as follows 
 $$ u_\e= \sum_{i: a_{i,\e}\in \partial \mathbb{S}^n_+} \frac{1}{K(a_{i,\e})^{(n-2)/4}}\d_{a_{i,\e},\l_{i,\e}} + \sum_{i: a_{i,\e}\in \mathbb{S}^n_+}  \frac{1}{K(a_{i,\e})^{(n-2)/4}} \d_{a_{i,\e},\l_{i,\e}} + v_\e,$$
 where  $\d_{a,\l}$ is the standard bubble defined in \eqref{eq:bble} and $||v_{\e}|| =  o_\e(1)$ where $||.||$ is defined by  \eqref{eq:norm}. \\
Furthermore, there hold:
\begin{enumerate}
  \item[(a)]
  Each interior concentration point  $a_{i,\e} \in \mathbb{S}^n_+$  converges to a critical point  $y_i$ of $K$ with  $\D K(y_i) < 0$,  $\l_{i,\e} d(a_{i,\e}, y_i) $ is uniformly bounded and    $y_i$ is an  \emph{isolated simple blow up point} for the sequence $(u_{\e})$. Moreover there exists  a dimensional constant $\kappa_1(n) > 0$ (see \eqref{kappa12} for the precise value) such that
  $$ \l_{i,\e} ^2 = - \kappa_1(n) \frac{\D K(y_i) }{ K(y_i)} \, \frac{1}{\e} (1+o_{\e}(1)).$$

  \item[(b)]
  Each boundary concentration point  $a_{i,\e}\in \partial \mathbb{S}^n_+$  converges to a critical point  $z_i$ of $K_1$ with   ${\partial_{\nu} K} (z_i) > 0$ and there exists a dimensional constant $\kappa_2(n)  > 0$  (see \eqref{kappa12} for the precise value) such that
 $$
 \l_{i,\e} \, = \, \kappa_2(n) \, \frac{\partial_{\nu}K(z_i)}{K(z_i)} \, \frac{1}{\e } (1+o_\e(1)).
 $$
  \item[(c)] For a boundary concentration point $a_{i,\e}$ converging to a critical point $z_i$ of $K_1$ there are two alternatives
  \begin{enumerate}
    \item[(i)] Either $ d(a_{i,\e},z_i)/ \e $ is uniformly bounded and $z_i$ is an  \emph{isolated simple blow up point} for the sequence $(u_{\e})$.
    \item[(ii)] Or for a subsequence $ d( a_{i,\e}, z_i) / \e \to \infty$,  as $\e \to 0$. In this case $z_i$ is a  \emph{non simple blow up point} for the sequence $(u_{\e})$ (in the sense that there exists at least another point  $a_{j,\e}$ which converges to $z_i$).
  \end{enumerate}
  \item[(d)] Let $z  \in \partial \mathbb{S}^n_+$ be  a non simple boundary blow up point  and  $\{a_{1,\e}, \cdots,a_{m,\e}\}$ be the set of points converging to $z$.  Then for every $i=1, \cdots,m$ we have that  $\e^{(2-n)/{n}} d(a_{i,\e},z_i)$ is uniformly bounded.\\
        Denoting
      \begin{equation}\label{bi}
        b_{i,\e}:= \ \kappa_3(n) \, \frac{(\partial_{\nu}K(z))^{\frac{n-2}{n}}}{(K(z))^{\frac{n-1}{n}}} \, \frac{1}{\e ^{\frac{n-2}{n}}}\, (a_{i,\e} -<a_{i,\e},z>z) \in  T_z \partial \mathbb{S}^n_+  \quad  \mbox{ for } 1\leq i\leq m,
      \end{equation}
     where $\kappa_3(n) > 0$ is a dimensional constant (see \eqref{kappa3} for the precise value),  we have that
       $(b_{1,\e},\cdots,b_{m,\e})$ converges to $(\ov{b}_1,\cdots,\ov{b}_m)$ which is a critical point of the following Kirchhoff-Routh type  function:
\begin{align} \label{dF2} & \mathcal{F}_{z,m}:  \mathbb{F}_m( T_z (\partial \mathbb{S}^n_+))\to \R ; \\
&   \mathcal{F}_{z,m}(\xi_1,\cdots,\xi_m):= \frac{1}{2}\sum_{i=1}^{m} D^2K_1(z)( \xi_i , \xi_i )  + \sum_{1 \leq i < j \leq m}\frac{1}{| \xi_i - \xi_j | ^{n-2}}, \nonumber \end{align}
where $ \mathbb{F}_m( T_z \partial \mathbb{S}^n_+):=  \{   (\xi_1,\cdots,\xi_m); \, \xi_i \neq \xi_j  \in T_z (\partial \mathbb{S}^n_+) \mbox{ for } i \neq j  \} $.
\end{enumerate}
\end{thm}

\begin{remark}
\begin{enumerate}
\item[(i)]
We point out that given  $z_1,\cdots,z_m \in \partial \mathbb{S}^n_+$  non degenerate critical points of $K_1:= K_{\lfloor \partial \mathbb{S}^n_+}$ satisfying  $(\partial_{\nu} K)(z_i) > 0$ for each $i=1,\cdots,m$ and  $y_{m+1},\cdots, y_{m+\ell}$  non degenerate critical points of $K$ with $\D K(y_i) < 0$ for each $i > m$,
then there exists a  sequence of solutions $u_{\e}$ of $(\mathcal{P}_\e)$  which converges weakly to $0$ and blows up at $z_1, \cdots, z_m, y_{m+1}, \cdots,y_{m+\ell}$ and all these blow ups are \emph{isolated simple}. See Theorem 1.1 in \cite{AB20b}.
  \item[(ii)]
Since the function $ \mathcal{F}_{z,m} $ does not have any critical point  if $z$ is a local maximum point, see Proposition 4.1 in \cite{AB20b}, it follows that $z$, in the situation of $(d)$ cannot be  local  maximum  of $K_1$. In other words  if a local maximum point $z$ is a blow up point then it has to be  isolated simple.
  \item[(iii)]
   The converse of  statement $(d)$ in the above theorem holds. Indeed  given $z \in \partial \mathbb{S}^n_+$ a critical point of $K_1$ having $\partial_{\nu} K(z) > 0$ and  $m \in \N$. Then   every non degenerate critical point of $\mathcal{F}_{z,m}$ gives rise to a solution of $(\mathcal{P}_{\e})$ building a cluster of $m$ concentration points around $z$. See Theorem 1.4  in \cite{AB20b}.
\end{enumerate}
\end{remark}

In the next theorem we describe the blow up picture when the sequence of solutions has a non zero weak limit. We call such a behavior a \emph{blow up phenomenon with residual mass}. Such a phenomenon does not occur in low dimension $n \leq 4$. We first describe the blow up scenario for half spheres of dimension greater than or equal to 7. Namely we prove:

\begin{thm}\label{th:t2} Let $n \geq 7$ and $0 < K \in C^3(\ov{\mathbb{S}^n_+})$   be a positive function  satisfying  the assumptions  $(H1), \, (H2)$ and let $(u_\e)$ be a sequence of energy bounded solutions of $(\mathcal{P}_\e)$ converging weakly but not strongly $u_\e \rightharpoonup \o \neq 0$. Then $\o$ is a solution of $(\mathcal{P})$ and $u_{\e}$  has to blow up and takes the following form
$$ u_\e \, = \,  \o \, + \,   \sum_{i: a_{i,\e}\in \partial \mathbb{S}^n_+} \frac{1}{K(a_{i,\e})^{(n-2)/4}}\d_{a_{i,\e},\l_{i,\e}} + \sum_{i: a_{i,\e}\in \mathbb{S}^n_+}  \frac{1}{K(a_{i,\e})^{(n-2)/4}} \d_{a_{i,\e},\l_{i,\e}} + v_\e,$$
where  $\d_{a,\l}$ is the standard bubble defined in \eqref{eq:bble} and $||v_{\e}|| = o_\e(1)$  where $||.||$ is defined by  \eqref{eq:norm}. \\
Furthermore the statements $(a)$,  $(b)$,  $(c)$ and  $(d) $   of Theorem \ref{th:t1} hold.
\end{thm}

\begin{remark}
  The converse of  Theorem \ref{th:t2} holds. See Theorems 1.2 and 1.7 in \cite{AB20b}.
\end{remark}
In the next theorem we single out the half sphere of dimension $5$. Since in this case, unlike the closed case, blow up with residual mass may occur but involves exclusively boundary blow up points. Namely we prove 

\begin{thm}\label{th:t3}
 Let $n = 5$ and  $0 < K \in C^3(\ov{\mathbb{S}^5_+})$   be a positive function  satisfying  the assumption  $(H1)$ and let $(u_\e)$ be a sequence of energy bounded solutions of $(\mathcal{P}_\e)$ converging weakly but non strongly  $u_\e \rightharpoonup \o\neq 0$. Then  $\o$ is a solution of $(\mathcal{P})$ and
 $u_\e$  blows up in the following form  
 $$ u_\e=  \o + \sum_{i: a_{i,\e}\in \partial \mathbb{S}^5_+} \frac{1}{K(a_{i,\e})^{3/4}}\d_{a_{i,\e},\l_{i,\e}}  + v_\e,$$
 where  $\d_{a,\l}$ is the standard bubble defined in \eqref{eq:bble} and $||v_{\e}|| = o_\e(1)$  where $||.||$ is defined by  \eqref{eq:norm}.\\
  Furthermore Statements  $(b)$, $(c)$ and  $(d)$ of Theorem \ref{th:t1} hold.
\end{thm}

\begin{remark}
  \begin{enumerate}
    \item[i)]
    The converse of  Theorem \ref{th:t3} holds. See Theorems 1.2 and 1.7 in \cite{AB20b}.
    \item[ii)] On $\mathbb{S}_+^6$  we constructed  in  \cite{AB20b} blowing up solutions with residual mass involving only boundary blow up points. However  the existence of  blowing up solutions with residual mass involving interior points remains an open problem even  for the Nirenberg problem on the six dimensional sphere $\mathbb{S}^6$.
  \end{enumerate}
\end{remark}

\begin{thm}\label{th:t4}
 Let $n = 5$ and  $0 < K \in C^3(\ov{\mathbb{S}^5_+})$   be a positive function  satisfying  the assumption  $(H1)$ or  $n \geq 7$ and $0 < K \in C^3(\ov{\mathbb{S}^n_+})$ satisfy the assumptions  $(H1)$ and $(H2)$. Let $(\o_k)$ be a sequence of energy bounded solutions of  $(\mathcal{P})$. Then $ | \o_k |_\infty$ is uniformly bounded (that is $(\o_k)$ cannot blow up).
 \end{thm}

\begin{remark}
We point out that the above results extend easily on compact riemannian  manifolds  with umbilic boundary. For more general manifolds there is a difficulty in the choice of a suitable \emph{modified bubble} when the boundary is not umbilic. This technical point will the subject of a forthcoming paper.
\end{remark}
Before closing this introduction we describe our strategy   of proof. To perform an asymptotic analysis of blowing up solutions $u_{\e}$ the usual blow up analysis techniques are based on \emph{ precise pointwise  $C^0$-estimates}  of $u_{\e}$ and the extensive  use of Pohozaev identities \cite{CL1, CL2, DHR, KMS09, yyli1, yyli2,  MM19, Schoen}. Our method of analysis is different. Indeed it  consists of   testing  the equation by vector fields which bring  the parameters of the concentration to their critical positions. These vector fields correspond to the leading term of the gradient of the Euler-Lagrange functional with respect to these parameters combined with what we call the \emph{barycentric vector field} which consists of pushing a cluster of nearby  blow up points to their common center of mass. We then   derive  balancing conditions to be satisfied by these parameters. Through a careful analysis of these balancing conditions,  we derive the information regarding the location of the concentration points and the speed  of the concentration. Eventually we derive the nature of the blow up point (i.e isolated  simple  or not). We believe that our method,  which avoids  the use of pointwise estimates and  Pohozaev identities   might be useful to deal with  non compact variational problems where \emph{non simple blow ups} occur like in the singular mean field equation with quantized singularities \cite{BT, KLin, WZ, WZb, DW}. Indeed the existence of non simple blow up points makes the task of establishing pointwise $C^0$-estimates a daunting one.
\bigskip


The remainder of the paper is organized as follows: in Section 2 we set up the variational framework of the problems $(\mathcal{P})$  and $(\mathcal{P}_{\e})$, introduce the neighborhood at infinity and its parametrization. Section 3 is devoted to the analysis of finite energy blowing up solutions in the  zero limit case. Indeed in  this section, after giving a precise estimates of the infinite dimensional part, we prove  various balancing conditions satisfied by  the parameters of the concentration and provide the proof of Theorem \ref{th:t1}. In Section 4 we deal with the non zero weak limit case and provide the proofs of Theorems \ref{th:t2} and \ref{th:t3}. While Section 5 is devoted to the proof of Theorem \ref{th:t4}. Finally we collect in the appendix some  useful estimates and technical lemmas needed  in the proof of various  statements in this paper.


\section{Parametrization of the neighborhood at infinity}

In this section we consider $\e \geq 0$ and  we will set up the general variational framework, recall the description of the lack of compactness that one derives from the concentration's compactness principle, introduce the \emph{neighborhood at infinity} and its parametrization.\\
The space of variation is  the Sobolev space  $H^1(\mathbb{S}^n_+)$  endowed with the norm
\begin{equation}\label{eq:norm}
||u||^2 := \int_{\mathbb{S}^n_+} |\n u|^2 \, + \, \frac{n(n-2)}{4} \int_{\mathbb{S}^n_+} u^2.
\end{equation}
Indeed, for $\e \geq 0$,  the  problem $(\mathcal{P}_\e)$ has a variational structure. Namely  its solutions are in one to one correspondence with the positive critical points of the functional
$$ I_\e(u) := \frac{1}{2}  \| u \|^2 - \frac{ 1 }{ p+1-\e}  \int_{ \mathbb{S}^n_+} K  | u |^{p+1-\e} , \quad u \in H^1(\mathbb{S}^n_+) \mbox{ with } p:= \frac{n+2}{n-2}.$$

For $a \in \ov{ \mathbb{S}^n_+}$ and $\l > 0$ we define the \emph{ standard bubble } to be
\begin{equation}\label{eq:bble}
\d_{a,\l}(x) := c_0 \frac{\l^{(n-2)/2}}{ (  \l^2 + 1 + (1-\l^2) \cos d(a,x))^{(n-2)/2}},
\end{equation}
where $d$ is the geodesic distance on  $\ov{\mathbb{S}^n_+}$ and $c_0:= (n(n-2))^{(n-2)/4}$ is chosen such that
$$
- \D \d_{a,\l} \, + \, \frac{n(n-2)}{4} \d_{a,\l} \, = \, \d_{a,\l}^{(n+2)/(n-2)} \quad \mbox{ in } \mathbb{S}^n_+.
$$

For $a \in \ov{ \mathbb{S}^n_+},$ we define \emph{projected bubble} $\varphi_{a,\l}$ to be the unique solution to
$$ - \D \varphi_{a,\l} \, + \, \frac{n(n-2)}{4} \varphi_{a,\l} \, = \, \d_{a,\l}^{(n+2)/(n-2)} \quad \mbox{ in }\,  \mathbb{S}^n_+; \quad \frac{\partial \varphi_{a,\l}}{\partial \nu} \, = 0 \mbox{ on } \partial \mathbb{S}^n_+.
$$
We point out that $\varphi_{a,\l} = \d_{a,\l}$ if $a \in  \partial \mathbb{S}^n_+$.
\medskip

Note that, the stereographic projection induces an isometry $\i : H^1(\mathbb{S}^n_+) \to D^{1,2}(\R^n_+)$ and the function $\i \d_{a,\l}$ is a solution of the Yamabe problem in $ \R^n_+$. Precisely,  for $a \in \partial \mathbb{S}^n_+$, using  the stereographic projection $\pi_{-a}$ (see \cite{BaBr} for the precise formulae), we have
$$ -\D (\i \d_{a,\l}) = (\i \d_{a,\l})^{(n+2)/(n-2)} \quad \mbox{ in } \R^n_+ \quad \mbox{ and } \quad \i \d_{a,\l}(x) := c_0 \frac{\l^{(n-2)/2}} { (1+\l^2 | x  |^2)^{(n-2)/2}}. $$
Note that, if we use $\pi_b$ (with $b\in \partial \mathbb{S}^n_+$ instead of $-a$) the image $\i \d_{a,\l}$ becomes $c_0 \mu^{(n-2)/2}/ (1+\mu^2 | x-\tilde{a} |^2)^{(n-2)/2}$ where the new variables $\tilde{a}$ and $\mu$ depend strongly on $a$ and $\l$ (see \cite{AB20b, BaBr} for the precise relations). In the sequel, we denote the function $\i\d_{a,\l}$ also by $\d_{a,\l}$.

Let $\o$ be a solution of $(\mathcal{P}_0)$  and let $N_0(\o)$ be the kernel of the associated quadratic form
defined by:
\begin{equation}\label{qw} Q_\o (h):= \| h \| ^2 - \frac{n+2}{n-2} \int_{\mathbb{S}^n_+} K \o ^{4/(n-2)} h^2 \qquad \mbox{ for } h \in H^1(\mathbb{S}^n_+). \end{equation}
 Let $m$ be the dimension of $N_0(\o)$ and  $(e_1,\cdots,e_m)$  be an orthonormal basis of $N_0(\o)$.  We set
$$
H_0(\o):= \, span(\o) \oplus span(e_1, \cdots,e_m).
$$
Following  M. Mayer,  see Lemma 3.6 and Proposition 3.7 in \cite{Mayer-Thesis}, we parameterize  a neighborhood of $\o$ by
\begin{align}
& u_{\a,\b} := \a (\o + \sum_{i=1}^m  \b_i e_i + h(\b))  \quad  \mbox{ with }  \label{ualphbeta} \\
&  h(\b) \perp H_0(\o); \, h(\beta) = O(||\b||^2) \mbox{ and } \| h(\b) \|_{C^2} \to  0, \nonumber \end{align}
where  the function  $u_{\a,\b} $ satisfies
\be \label{ualphbeta2} \langle  \n I_0(u_{\a,\b}), h \rangle := \langle u_{\a,\b},h\rangle  - \int K u_{\a,\b}^{p} h = 0  \quad \mbox{ for each } h\in H_0(\o)^{\perp} .\ee

Next for $\o$ a solution of $(\mathcal{P})$ whose Kernel is of dimension $m$,  $q,\, \ell \in \N_0$ and $\tau$ a small positive constant, we define the so called \emph{neighborhood at infinity} $V(\o, q,\ell,\tau)$ as follows:

\begin{align*}
{V}(\o, q, & \ell, \tau)  :  =   \Big \{  u \in H^1(\mathbb{S}^n_+):  \,
 \exists \,  \l_1, \cdots, \l_N > {\tau^{-1}} \mbox{ with } \e \ln \l_i \leq \tau ; \,  \exists \, a_1, \cdots,a_{N} \in \ov{ \mathbb{S}^n_+}, \\
 &  \mbox{ with }  \l_i d_i  < \tau,   \forall \,  i \leq q, \,   \mbox{and } \l_i d_i > \tau^{-1} \, \,  \forall \,   i > q,  \, \, \,    \e_{ij} < \tau ; \, \,
   \exists \, \a_0 \in(1-\tau,1+\tau) ;\\
   &  \b \in \R^m   {  \mbox{ with }\| \b \| \leq c\, \tau }\,  \mbox{ such that } \| u -  \sum_{i = 1}^N K(a_i)^{\frac{2-n}{4}} \varphi_{a_i,\l_i} - u_{\a_0,\b} \| < \tau  \Big \}
\end{align*}

where $N:= q+\ell$,  $d_i := d(a_i,\partial \mathbb{S}^n_+)$ and
\begin{equation}\label{eppe}
\e_{ij}:= \Big( \frac{\l_i}{\l_j} + \frac{\l_j}{\l_i} + \frac{1}{2} \l_i \l_j (1 - \cos  d (a_i,a_j))\Big)^{(2-n)/2}.
\end{equation}

Following A. Bahri and J-M. Coron \cite{BCd}  we consider for  $u \in V(\o, q,\ell,\tau)$  and $N= q + \ell$ the following minimization problem

\begin{equation}\label{eq:min}
  \min_{  \a_i > 0 ; \, \, \b\in \R^m ; \, \,   \l_i > 0, \, \, a_i \in \partial \mathbb{S}^n_+, \forall i \leq q; \, \,  a_i \in \mathbb{S}^n_+, \forall  i > q  }     \Big\|  u  \, - \, \sum_{i=1}^{N}\a_i \varphi_{a_i,\l_i} - u_{\a_0,\b} \Big\|  .
\end{equation}
We then have the following proposition whose proof is identical, up to minor modification to  the one  of  Proposition 7 in \cite{BCd}

\begin{pro}\label{p:min}
For any $q,\ell  \in \N_0$ there exists $\tau_0 > 0$ such that if $\tau < \tau_0$ and $u \in V( \o, q, \ell, \tau)$ the minimization problem \eqref{eq:min}  has, up to permutation of the indices, a unique solution.
\end{pro}
Hence it follows from Proposition \ref{p:min} that every $u \in V(\o,  q,\ell, \tau)$ can be written in a unique way as
\begin{align}
& u \, = \, \sum_{i=1}^{q} \a_i \d_{a_i,\l_i} \, + \,  \sum_{i= q+1}^{N} \a_i \varphi_{a_i,\l_i} \, + \, u_{\a_{0},\b} \, +  \, v, \qquad \mbox{ where }\label{F1}\\
& a_i \in \partial \mathbb{S}^n_+, \, i \leq q \mbox{ and } a_i \in \mathbb{S}^n_+ \, \, \mbox{ with } \, \l_i d_i  > \tau^{-1}\,  \, ,\, \,   i \geq q+1,\nonumber \\
& \a_i^{4/(n-2)} K(a_i) = 1 + o(1) \, \,  \forall \, i\geq 1,  \quad \a_0 = 1 + o(1),  \label{alphai0}
\end{align}
and  $v \in H^1(\mathbb{S}^n_+)$ satisfying
\be (V_0) \qquad  \| v\| <  \tau, \quad <v, \psi> = 0, \mbox{ for } \psi \in  E_{\o, a,\l}^\perp \mbox{ where } \label{eq:V0} \ee
\be E_{\o, a,\l}:=\mbox{span} \{  \d_i, \frac{\partial \d_i}{\partial \l_i}, \frac{\partial \d_i}{\partial a_i},  \varphi_j, \frac{\partial \varphi_j}{\partial \l_j}, \frac{\partial \varphi_j}{\partial a_j} , \, u_{\a_0,\b} , \frac{\partial u_{\a_0,\b}}{\partial \b_k} ;   1 \leq i\leq q; \, q  < j \leq N ; \, k\leq m\} ,   \label{Ewal} \ee
where $ \d_i := \d_{a_i,\l_i}$ and $\varphi_i:= \varphi_{a_i,\l_i}$.


Finally by combining the analysis of the Palais-Smale sequences, performed  in Lemma 3.3 of \cite{PT95}, whose proof goes along  with  the concentration compactness arguments  developed  in \cite{Struwe, Lions, BrC, DHR}, with the above parametrization of the neighborhood at infinity we derive that:

\begin{pro}\label{lambdaepsilon}
Let $ u_\e$  be an energy bounded solution of $(\mathcal{P}_\e)$ which  converges weakly to $0$. Then there exist $q$ and $\ell$ such that $u_\e$ can be written as
$$ u_\e:=  \sum_{i\leq q} \a_{i,\e} \d_{a_{i,\e},\l_{i,\e}} + \sum_{i  =q+1}^{q+\ell} \a _{i,\e} \varphi _{a_{i,\e}, \l_{i,\e}} + v_\e \in V(q,\ell ,\tau ):= V(0, q,\ell ,\tau )$$
 with  $ v_\e \in E_{a_\e ,\l_\e }^\perp:= E_{0,a_\e ,\l_\e }^\perp$.
\end{pro}



\section{The  zero weak limit case}\label{ss3}
In this section we deal with the case where the blowing up sequence $u_{\e}$ of $(\mathcal{P}_{\e})$  has a bounded energy and converges weakly to zero.   We notice that, in the sequel, for sake of simplicity of the presentation, we will cancel the index $\e$ from the points $a_{i,\e}$ and the speeds $\l_{i,\e}$.  \\
We remark that, in the first part  of this section, we consider $\e \geq 0$. That is the sequence $(u_\e)$ (in the case of $\e=0$) will be assumed to be a sequence $(\o_k)$ of solutions of $(\mathcal{P})$. However, in Proposition \ref{cp} and Lemma \ref{nonsimple1}, we will assume that $\e > 0$. 

\subsection{ Estimates of the infinite dimensional part}

In this subsection we deal with  $v_\e$ the infinite dimensional part of the blowing up solutions $u_{\e}$.  In the next lemma we provide an accurate estimate of this part, which shows that it does not have any contribution to  the blowing up  phenomenon. Namely we prove

\begin{pro}\label{eovv} Let  $v_\e$  be defined as in Proposition \ref{lambdaepsilon}. Then it satisfies : $$ \| v_\e \| \leq c \,  R(\e, a, \l)    \qquad \mbox{ where }$$
$$ R(\e,a,\l):= \e +  \sum_{i=1}^N \frac{| \n K(a_i) | }{\l _i} + \frac{1}{\l _i^2} + \begin{cases}
 \sum \e _{ij }^{\frac {n+2}{2(n-2)}}(\ln \e _{ij}^{-1})^{\frac{n+2}{2n}} + \sum_{i > q }\frac{\ln(\l_i d_i)}{(\l_i d_i)^{\frac{n+2}{2}}} & \mbox{ if } n \geq 6,\\
  \sum \e _{ij }(\ln \e _{ij}^{-1})^{{3}/{5}} + \sum_{i > q} \frac{1}{(\l_i d_i)^{3}} & \mbox{ if } n = 5. \end{cases}
$$
\end{pro}
\begin{pf}
Observe that, for each $u, h \in H^1(\mathbb{S}_+^n)$, it holds that
\be\label{www1} \langle \n I_\e (u), h \rangle = \langle u , h \rangle - \int K | u |^{p-\e -1} u h.\ee
Hence taking $u=u_\e$ and $h=v_\e$ and using the fact that $v_\e \in E_{a,\l}^\perp$, we derive that
\be\label{ff6} \| v_\e \|^2 = \int K u_\e ^{p-\e} v_\e = \int K \ov{u}_\e ^{p-\e} v_\e + p \int K \ov{u} _\e ^{p-1-\e} v_\e^2  + o(\| v_\e\|^2 ) \ee  where $\ov{u}_\e := u_\e - v_\e.$
Since $ \a_i^{4/(n-2)} K(a_i) = 1 + o(1) $ for each $i$, easy computations imply that
$$  \| v_\e \| ^2  -  p \int K \ov{u} _\e ^{p-1-\e} v_\e^2 =  \| v_\e \| ^2 - p \sum \int \d_i ^{p-1} (v_\e) ^2 +  o(\| v_\e\|^2 ) := Q_{a,\l}(v_\e) +  o(\| v_\e\|^2 ). $$
We remark that $Q_{a,\l}$ is a positive definite quadratic form in the space $E_{a,\l}^\perp$ (see Proposition 3.1 of \cite{B1}). Furthermore, the linear form is computed in Eq. (23) of \cite{AB20b} and we have
\be\label{www2} \int_{\mathbb{S}^n_+} K \ov{u}_\e ^{p-\e} v_\e = O \Big( R(\e, a, \l)  \| v_\e  \|\Big) . \ee

Combining the previous estimates we obtain
\begin{align*} Q_{a,\l}(v_\e) & = o(  \| v_\e\|^2) + O\Big( \| v_\e\| \Big( R(\e,a,\l)  \Big)\Big). \end{align*}
The estimate of $v_\e $ follows from the fact that $Q_{a,\l}$ is a positive definite quadratic form.
\end{pf}

\subsection{Balancing conditions for the parameters of the concentration}

We next  provide various balancing conditions that have to be satisfied by the parameters of concentration.  The following propositions are quoted from \cite{AB20b} (Propositions 3.1 and 3.2) by using the fact that $\n I_\e(u_\e)= 0$ since $u_\e$ is a solution of $(\mathcal{P}_\e)$.
Our first such a condition concerns the gluing parameter $\a_i$ for $i=1, \cdots,N$. Namely we have:
\begin{pro}\label{p:35} \cite{AB20b}
For each $1\leq i \leq N$, it holds:
\begin{align}\label{alpha1} & | 1- \l_i^{-\e(n-2)/{2}} \a _i^{{4}/{n-2}}K(a_i) | = O( R_{\a_i}) \qquad \mbox{ where }\\
 & R_{\a_i}:= \e +  \frac{ 1 }{\l_i^2 }+ \sum_{j\neq i} \e_{ij}   +  R(\e,a,\l)^2   + \begin{cases}   { | \n K(a_i) | } / {\l_i } \quad \mbox{if }i \leq q,   \\
  {1} / {(\l_i d_i)^{n-2}}   \quad \mbox{if }i \geq q+1 \end{cases}  . \end{align}
\end{pro}

Our next balancing condition concerns the rate of the concentration $\l_i$ for $1 \leq i \leq N$. Using Propositions 3.1, 3.2 of \cite{AB20b} and Proposition \ref{p:35} above, we get  :

\begin{pro}\label{p:33}  \cite{AB20b}
  For $\e$ small enough, the following equations  hold:
\be\label{(Ei)} (E_i) \begin{cases} \displaystyle{ - \frac{c_2} {2}\sum _{i \neq j\leq q} \a _j \l _i\frac{\partial \e _{ij}}{\partial \l _i} }  & -  \displaystyle{ \, \,  \frac{\a_i }{K(a_i)} \Big[ \frac{c_3}{\l _i}\frac{\partial K}{\partial \nu}(a_i) -  c_5 K(a_i) \e\Big] } \\  & =  \displaystyle{ O ( \frac{1}{\l_i ^2} + \sum_{j > q} \e_{ij} + R_{1,\l}   ) } ; \quad    \mbox{if } i \leq q \\
\displaystyle{  - {c_2} \sum _{j\neq i} \a _j \l _i\frac{\partial \e _{ij}}{\partial \l _i} } & + \displaystyle{\, \, \,  {c}_2  \frac{n-2}{2}   \sum_{j=q+1}^N\a_j  \frac{H(a_i,a_j)}{(\l_i\l_j)^{\frac{n-2}{2}}}  }\\
& + \, \, \displaystyle{  {\a_i} \Big(c_4 \frac{\Delta K(a_i)}{\l_i^2 K(a_i)} + 2c_5 \e\Big)  =   O (  R_{2,\l}  ) } ; \, \quad    \mbox{if } i > q \end{cases} \ee

 where $ R_{1,\l} :=  \sum_{j\neq k}  \e _{kj }^{\frac{n}{n-2}}\ln (\e _{kj }^{-1})  \ + R^2(\e,a,\l), $  $ R_{2,\l}:=  R_{1,\l} + \sum_{k > q} \ln(\l_k d_k) / (\l_k d_k)^n$ and
 \begin{align*}
  &  c_2:= c_0^{p+1}\int_{\R^n}\frac{1}{(1+|x|^2)^{(n+2)/2}}dx ; \, \, \,   c_3 :={(n-2)}\, c_0^{\frac{2n}{n-2}}\int_{\R_+^n}\frac{x_n(|x|^2-1)}{(1+|x|^2)^{n+1}}dx, \\
&  c_4:= \frac{n-2}{2n} c_0^{p+1} \int_{\R^n} \frac{ | x |^2( |x|^2 - 1 )}{(1+| x |^2)^{n+1}} ; \, \, \, \, \,   c_5 := \frac{(n-2)^2}{4} c_0^{p+1}\int_{\R_+^n}\frac{(|x|^2-1) \ln(1+|x|^2)}{(1+|x|^2)^{n+1}}dx .\end{align*}
\end{pro}

Our last balancing condition concerns the point of concentration $a_i$ for $i=1, \cdots,N$. Namely we prove
  \begin{pro} \label{p:38}  \cite{AB20b}
   For $\e$ small enough, there hold:
\be \label{Fib}  (F_i) \quad  \begin{cases}
\displaystyle{ - \frac{c_2}{2} \sum_{i \neq j\leq q} \a _j  \frac{1}{\l _i}\frac{\partial\e _{ij}}{\partial a_i } -  \frac{ \a _i} {K(a_i)}\frac{c_6}{\l _i}{\n K_1}(a_i) = O (  \frac{1}{\l_i ^2} + \sum_{j > q} \e_{ij} + R_{a_i}   ) } ; \, \,  \mbox{for } i \leq q \\
\displaystyle{  \frac{|\nabla K(a_i)|}{\l_i} \leq c \Big(\frac{1}{\l_i^3}+ \frac{1}{(\l_id_i)^{n-2}} + \sum_{j\neq i}\e_{ij} + R^2(\e,a,\l) \Big) } ; \, \,  \mbox{for } i > q \end{cases}\ee 
  where
  $$ R_{a_i} := R_{1,\l} + \sum_{j\leq q; j \neq i} \e _{ij}^{\frac{n+1}{n-2}}\l _j d(a_i , a_j )  \quad , \quad  c_6=\frac{n-2}{n}c_0^{\frac{2n}{n-2}}\int_{\R^n_+}\frac{| x |^2}{(1+|x|^2)^{n+1}}dx.$$
\end{pro}


\subsection{ Refined blow up analysis}

 In this section we analyze the nature of the blow ups. Namely we will prove that while  interior blow ups are \emph{isolated simple} and their  rates of concentration are comparable, there might be \emph{non simple boundary blow up points}. The existence of such \emph{non simple blow up points}  at a boundary blow up point $z \in \partial \mathbb{S}^n_+$ is in one to one correspondence with the  existence of critical points of some  Kirchhoff-Routh type  Hamiltonian $\mathcal{F}_z$. \\
 To perform our analysis we make use of the following notation.  We denote by \be\label{mu1} I_b:= \{i: a_i \in \partial \mathbb{S}^n_+\} ; \quad I_{in}:= \{i: a_i \in \mathbb{S}^n_+\} ; \quad \mu_i:= \begin{cases} \l_i \mbox{ if } i\in I_b\\
                       \l_i^2 \mbox{ if } i\in I_{in}\end{cases} . \ee
 Note that, in our case, $\#I_b = q$ and $\# I_{in} = \ell$. Furthermore, we order the $\mu_i$'s as \be \label{mu2} \mu_1 \leq \cdots \leq \mu_N.\ee Let
 \be \label{I} I' := \{ i : \lim_{\e \to 0} \mu_i / \mu_1 = \infty\} \quad ; \quad I:= \{1,\cdots,N\} \setminus I' \quad \mbox{ where } N:= q+\ell.\ee
 We notice that, for $i,j \in I$, the    speeds  of concentration  $\mu_i$ and $\mu_j$ are of the same order (i.e. the ratio is bounded above and below) and for each $k \notin I$ and each $ i \in I$, we have  that the ratio $\mu_k/\mu_i$ goes to $\infty$ as $\e\to 0$. Finally, for a critical point  $z$ of $K_1$ (resp. a critical point $y$ of $K$) we denote by
 \be\label{Bzy} B_z:= \{ i \in I \cap I_b: \lim_{\e \to 0} a_i = z\} \quad ; \quad  B_y:= \{ i \in I \cap I_{in}: \lim_{\e \to 0} a_i = y\} .\ee

Before studying the nature of the blow up we need the following Lemmata.

 \begin{lem}\label{sum} For each solution $u_\e$ in $V(q,\ell,\tau)$, it holds
 $$\sum_{k\neq i} \e_{ik} + \sum_{k\in I_{in}} \Big( \frac{ | \n K(a_k) | }{ \l_k} + \frac{1}{(\l_k d_k)^{n-2}}\Big)  + \e \leq c  \frac{1}{\mu_1}.$$
 Furthermore $R_{1,\l}$ and $R_{2,\l}$ defined in Proposition \ref{p:33} satisfy: $R_{k,\l} =  o\Big(\frac{1}{\mu_1^{(n-1)/(n-2)}}\Big)$.
 \end{lem}
 \begin{pf}
 We first notice that
 \begin{align}\label{derepsilon}
& - \l_i \frac{\partial \e_{ij}}{\partial \l_i} - \l_j \frac{\partial \e_{ij}}{\partial \l_j
} \geq 0 \, \, \mbox{ for each } i\neq j  \, \,  \mbox{ and } \, \,   - \l_i \frac{\partial \e_{ij}}{\partial \l_i} \geq c \e_{ij} \quad \mbox{ if } \l_i \geq c' \l_j\\
& -  \l_i \frac{\partial \e_{ij}}{\partial \l_i} \geq c \e_{ij} \quad \mbox{ for each } i\in I_{in} , \, \, j \in I_b.\end{align}
Hence, summing $2^i (E_i)$ (defined in \eqref{(Ei)}) for $i\in I_{in}$ and  for $i \in I_b$ respectively, we obtain
\begin{align*}
 & \sum_{i\in I_{in}; k\neq i} \e_{ik} +  \sum _{i\in I_{in}} \frac{1}{(\l_i d_i)^{n-2}} +  \e  \leq c   \Big( R_{2,\l} + \frac{1}{\mu_1}\Big) \\
  & \sum_{i, k \in I_b; k\neq i} \e_{ik} +  \e  \leq c   \Big( R_{1,\l}  +R_{2,\l}+ \frac{1}{\mu_1}\Big).\end{align*}
The result follows by using $(F_i)$ for $i \in I_{in}$ and the definitions of $R_{1,\l}$ and $R_{2,\l}$.
\end{pf}

 \begin{lem}\label{56}
 For $i\in I_{in}$, we define  \be\label{Gi}\G_i:= \l_i^2 \sum_{k \neq i} \e_{ik} + \frac{ H(a_i,a_i)}{\l_i ^{n-4}} + \frac{| \n K(a_i) | /\l_i }{\sum \e_{ki} + 1/\l_i^2}.\ee
 Let $D_1':=\{ i \in I_{in}: \lim \G_i =\infty\}$ and $D_1:= I_{in} \setminus D_1'$. \\
$(i)$  For each $i \in D_1$, there exists a critical point  $y_i$ of $K$ such that $\l_i | a_i - y_i | \leq C$.\\
$(ii)$ If there exists an index $i \in D_1'$, then it holds that
$$ \sum_{j\in I_{in}; j \geq i} \Big( \frac{ | \n K(a_i) | }{\l_i} + \sum_{k\neq j} \e_{kj} + \e  + \frac{1}{(\l_j d_j)^{n-2}} \Big) = o\Big(\frac{1}{\mu_1^{(n-1)/(n-2)}}\Big). $$
 \end{lem}
 \begin{pf}
 Let $i\in D_1$, it follows that $\G_i$ is bounded which implies  that $| \n K(a_i) | \leq C/\l_i$ and therefore  the first assertion follows. Concerning the second one, let $i\in D_1'$, summing $2^j (E_j)+(F_j)/m$ for $j \geq i$ and $j\in I_{in}$ where $m$ is a small constant, it holds that
 $$ \sum_{j\in I_{in}; j \geq i} \Big( \frac{ | \n K(a_j) | }{\l_j} + \sum_{k\neq j} \e_{kj} + \e  + \frac{H(a_j,a_j)}{\l_j ^{n-2}} \Big) = O\Big( \frac{1}{\l_i^2} + R_{2,\l} \Big).$$
 Since $i\in D_1'$, it follows that $1/\l_i ^2$ is small with respect to the left hand side. Thus the proof follows from the estimate of $R_{2,\l}$ (see Lemma \ref{sum}).
 \end{pf}

 \begin{lem}\label{55}
 $(i)$ For each $i \in I_{in}$ and each $j \in I_b$, it holds that : $\e_{ij}= o(1/\mu_1^{(n-1)/(n-2)})$.\\
 $(ii)$ For each $i, j \in I_{in}$, it holds that : $\e_{ij}= o(1/\mu_1)$.
 \end{lem}
 \begin{pf}
Claim $(i)$ follows from Lemma \ref{56}. It remains to prove the second claim. Let $i,j \in I_{in}$. Observe that, if $i$ or $j$ belongs to $D_1'$, then the result follows from Lemma \ref{56}. In the other case, that is $i,j \in D_1$, using again Lemma \ref{56}, there exist critical points $y_i$ and $y_j$ such that $\l_k d( a_k, y_k ) \leq c $ for $k=i,j$. Two cases may occur:
\begin{enumerate}
\item  [(a)] Either $y_i \neq y_j$, and in this case we get $d( a_i, a_j ) \geq c $ and therefore the result follows easily,
\item [(b)] or $y_i = y_j$. Since we have $\l_k  d( a_k, y_k ) \leq c $ and $\e_{ij}$ is small, it follows that $\l_i / \l_j \to 0$ or $\infty$. Taking $\l_i \leq \l_j$ and using the fact that $\G_j$ is bounded, it holds
 $$ \e_{ij} \leq \frac{c}{  \l_j ^2}  = c \frac{ \l_i ^2}{  \l_j ^2} \frac{1}{  \l_i ^2} = o \big( \frac{1}{  \l_i ^2} \big) = o \big( \frac{1}{  \mu_1} \big).$$
 \end{enumerate}
 Hence  the proof is completed.
 \end{pf}

\subsubsection{Ruling out bubble towers}

In this section we prove that the rate of concentration of boundary concentration points are comparable and the rate of concentration of interior points are also comparable. This fact rules out the phenomenon of \emph{bubble towers}.\\
We start with a preliminary lemma:
  \begin{lem}\label{57}
$(1)$ Assume that there exist $i_0 \neq j_0 \in I_b$  such that $\lim \l_{i_0} \e_{i_0j_0} = \infty$. Then it holds that
$$\e +  \sum_{k\in I_b; k\neq i_0} \e_{i_0k} + \frac{1}{\l_{i_0}} = o\Big( \frac{1}{\mu_1^{(n-1)/(n-2)}}\Big).$$
 $(2)$ For each $i,j\in I_b$ with $j\neq i$, it holds
 \begin{align*}
 &(a)\quad  \frac{1}{\l_i} | \frac{\partial \e_{ij}}{\partial a_i} | \leq {c}/{\l_i ^{(n-1)/(n-2)} }\quad \mbox{ if } \quad c' \leq \frac{\l_i}{\l_j} \leq c \, \, \mbox{ and } \, \,  \l_i \e_{ij} \leq c,\\
 & (b) \quad  \frac{1}{\l_i} | \frac{\partial \e_{ij}}{\partial a_i} | = o\big( 1/{\mu_1^{(n-1)/(n-2)}} \big) \quad \mbox{ in the other cases}.
 \end{align*}
   \end{lem}
  \begin{pf}
Summing $2^i (E_i)$ for $i\in I_b$ and $i\geq i_0$, it follows that
$$ \sum_{i\in I_b; i\geq i_0} \sum_{  k\in I_b; k\neq i} \e_{ik} + \e + O\big( \frac{1}{\l_{i_0}} \big) =O\big(\sum_{k\in I_b; j \in I_{in}} \e_{kj} + R_{1,\l} \big).$$
  Hence assertion $(1)$ follows from Lemmas \ref{sum}, \ref{55} and the fact that $1/\l_{i_0}$ is small with respect to $\e_{i_0j_0}$ which exists in the left hand side.\\
  Now we will focus on the second assertion. Observe that, in general, it holds
  \be\label{c2} \frac{1}{\l_i} | \frac{\partial \e_{ij}}{\partial a_i} | \leq c \,  \l_j d( a_j , a_i ) \e_{ij}^{n/(n-2)} \leq c \, \sqrt{{\l_j}/{\l_i}} \e_{ij}^{(n-1)/(n-2)} \leq c\, \,  \e_{ij} .\ee
  Observe that, Claim $(a)$ follows from the second inequality of \eqref{c2}. Concerning Claim $(b)$, three cases may occur:
  \begin{itemize}
\item if $\l_i \e_{ij} \to \infty$ or $\l_j \e_{ij} \to \infty$, then $(b)$ follows from \eqref{c2} and Claim $(1)$.
\item if $(\l_i  +\l_j)\e_{ij} \leq c$ and $\l_i /\l_j \to \infty$, then $(b)$ follows from the second inequality of \eqref{c2}.
\item if $(\l_i  +\l_j)\e_{ij} \leq c$ and $\l_i /\l_j \to 0$, then it holds
$$ \sqrt{{\l_j}/{\l_i}} \e_{ij}^{(n-1)/(n-2)} \leq  (\l_j\e_{ij})^{\frac{n-1}{n-2}} (\l_i/\l_j)^{\frac{n-1}{n-2} - \frac{1}{2}}(1/\l_i)^{\frac{n-1}{n-2}} = o({1}/{\mu_1^{(n-1)/(n-2)}}).$$
\end{itemize}
Hence the proof of $(b)$ is completed.
  \end{pf}

 Now, using the above lemmas, Equations $(E_i)$ and $(F_i)$ can be improved and we get
 \be\label{Ei'} (E_i') \, \begin{cases}
 \displaystyle{  2c_5 \e + c_4 \frac{ \D K(a_i) }{\l_i ^2 K(a_i)} = o\big( \frac{1}{\mu_1} \big) } & \mbox{ for } i\in I_{in}, \\
 \displaystyle{ - \frac{c_2} {2}\sum _{ i \neq  j\in I_b} \a _j \l _i\frac{\partial \e _{ij}}{\partial \l _i}  -  \frac{\a_i }{K(a_i)} \Big[ \frac{c_3}{\l _i}\frac{\partial K}{\partial \nu}(a_i) - c_5 K(a_i) \e\Big] = o\Big( \frac{1}{\mu_1^{\frac{n-1}{n-2}}} \Big)}  & \mbox{ for } i\in I_{b} , \end{cases} \ee
 \be \label{Fi'} (F_i')  : \quad   - \frac{c_2}{2} \sum_{j\in I_b; j\neq i} \a _j  \frac{1}{\l _i}\frac{\partial\e _{ij}}{\partial a_i } -  \frac{ \a _i} {K(a_i)}\frac{c_6}{\l _i}{\n K_1}(a_i) = o \Big( \frac{1}{\mu_1^{(n-1)/(n-2)} } \Big)  \quad \mbox{ for } i\in I_{b} . \ee

Next  we rule out bubble towers by proving that all  concentration's rates $\mu_i$ are comparable. Namely we prove
\begin{pro}\label{comparable}
All the $\mu_i$'s are comparable. That is
$$ I = I_b \cup I_{in} .$$
\end{pro}
\begin{pf} Recall that $N= q + \ell$.
Arguing by contradiction, we assume that $ N \notin I$. Thus we get that $\lim \mu_{ N }/\mu_1 = \infty$. We claim that\\
{\bf Claim 1:}  $ \e = o(1/\mu_1)$. In fact two cases may occur: $(a)$ Either $ N \in I_{in}$, and in this case the claim follows easily by using $(E_N')$ or $(b)$ $ N \in I_b$, in this case, $ (E_N')$ implies that
$$  \sum_{j\in I_b; j\neq N} \e_{jN} + \e  + O(\frac{1}{\l_{ N }}) = o\big( \frac{1}{\mu_1^{(n-1)/(n-2)}} \big) $$ and the claim follows in this case also. Hence the proof of Claim 1.\\
{\bf Claim 2:} $ I\cap I_{in } = \emptyset$.  Assume that there exists $j \in I\cap I_{in}$. It follows that $\mu_1$ and $\l_j^2$ are of the same order. Hence using Claim 1 and $(E_j')$ we derive that $\D K(a_j) =o(1)$. On the other hand,  using Lemma \ref{sum}, we derive that $a_j$ has to converge to a critical point of $K$ which leads to a contradiction. Thus Claim 2 follows. \\
{\bf Claim 3:} For each $i \in I$, we have: $a_i $ converges to a critical point $z_i$ of $K_1$ satisfying $\partial K/ \partial \nu (z_i) > 0$.
Using Claim 2, it follows that $I\subset I_b$. Let $i\in I$, using $(F_i')$ and Lemma \ref{57}, we get $| \n K_1 (a_i) | /\l_i = o(1/\mu_1)$ which implies that $| \n K_1 (a_i) |=o(1)$ since $\l_i$ and $\mu_1$ are of the same order. Thus $a_i$ has to converge to a critical point $z_i$ of $K_1$. It remains to prove that  ${\partial_{\nu} K}(z_i) > 0$.  First, we remark that
\be\label{c3} \e_{kj} = o(1/\mu_1) \quad \mbox{ for each } k\in I \mbox{  and } j \in I_b \setminus  I\, \,  (\mbox{that is } \l_j/\l_k \to \infty). \ee
 In fact, \eqref{c3} follows from Lemma \ref{57} if $\l_j \e_{kj} \to \infty$. In the other case, that is $\l_j \e_{kj} \leq c$, it holds that $\e_{kj} \leq c/ \l_j =o(1/\l_k)$ and \eqref{c3} follows in this case also. \\
 Secondly, using \eqref{c3}, Claim 1  and $(E_i')$, we get
$$ \sum_{ k\in I; k\neq i} \a_k (c+o(1))\e_{ik} - c\, \a_i \,  \frac{\partial K}{\partial \nu}(z_i)\frac{1}{\l_i} = o\big(\frac{1}{\mu_1}\big)$$ which implies that ${\partial_{\nu} K}(z_i)$ has to be positive. Thereby Claim 3 is completed.\\
{\bf Claim 4:} For each critical point $z$ of $K_1$, it holds $\# B_z \neq 1$ where $B_z:= \{ i\in I\cap I_b : \lim a_i =z\}$. \\
Assume that there exists $z$ such that $\# B_z =1$. Let $ B_z = \{i\}$.  It follows that $\e_{ij}= o(1/\mu_1)$ for each $j\in I_b$ with $j\neq i$ (in fact, for $j \in I$, we derive that $| a_i - a_j| \geq c$ and for $ j\in I_b \setminus I$, it follows from \eqref{c3}). Thus $(E_i')$ leads to a contradiction (by using Claim 1) and therefore Claim 4 follows.

In the sequel, let $z$ be such that $\# B_z \geq 2$ and let $i_0,i_1\in B_z$ be such that $| a_{i_0 }- a_{i_1} | := \min \{ | a_i-a_j|, i,j\in B_z \mbox{ with } i\neq j\}$. We introduce the following sets: $$ A_z':= \{ j \in B_z: \lim | a_j - a_{i_0} | / | a_{i_0 }- a_{i_1} | = \infty\} \quad \mbox{ and } \quad A_z := B_z \setminus A_z'.$$
We remark that $A_z$ contains at least $i_0$ and $i_1$.\\
{\bf Claim 5:} Let $z$ be such that $\#B_z \geq 2$, it holds that $(i)$: $\sum_{k\neq i; k\in A_z} \e_{ik} = (c +o(1))/\l_i$ for each $i \in A_z$. Furthermore, $(ii)$: for each $i\neq j \in A_z$, it holds $\l_i ^{(n-3)/(n-2)} d(a_i,a_j) $ is bounded above and below.\\
We notice that, for each $i\neq j \in A_z$, it holds that $d(a_i,a_j)$ and $d(a_{i_0}, a_{i_1})$ are of the same order. Furthermore, $\l_i$ and $\l_j$ are of the same order. Hence, $\e_{ij} = (\l_i \l_j d(a_i,a_j)^2)^{(2-n)/2} (c+o(1))$ which implies that
all the $\e_{ij}$, for each $i\neq j \in A_z$, are of the same order. Thus $(ii)$ follows immediately from $(i)$. Concerning $(i)$, it follows from Claim 1, $(E_i')$ and the fact that $\e_{ij}= o(1/\mu_1)$ for each $i\in A_z$ and $j\notin A_z$ (the last information is immediately for $j \in I\setminus B_z$. It follows from \eqref{c3} for $j \notin I$.  For $j \in A_z'$, we have $d(a_i,a_j) / d(a_{i_0},a_{i_1}) \to \infty$ and the $\l_k$'s are of the same order which implies that $\e_{ij}=o(\e_{i_0 i_1})$). Hence Claim 5 is complete.

\medskip

To conclude the proof of Proposition \ref{comparable}, we need to introduce the barycenter of the points $a_i$'s for $i \in A_z$. Let $b \in \R^{n+1}$ be such that $\sum_{i\in A_z} (b- a_i) = 0$ and we define  $\ov{ a } := \frac{ b }{ | b | }$. This point is the barycenter of the points $a_i$'s and it satisfies
\be\label{rty1} \ov{ a } \in \partial \mathbb{S}^n_+ \quad ; \quad a_i - \langle a_i,\ov{a} \rangle \ov{a} \in T_{\ov{a} } ( \partial \mathbb{S}^n_+) \, \, \forall \, \, i \in A_z \quad \mbox{ and } \quad  \sum_{i\in A_z} a_i - \langle a_i,\ov{a} \rangle \ov{a} = 0 .\ee
In addition it is easy to get that $\l_i | \ov{a} -\langle a_i,\ov{a}\rangle a_i | \leq c \l_i d(a_i, \ov{a} ) \leq c \l_i ^{1/(n-2)}$ (by using Claim 5).
Now, multiplying $(F_i')$ by $\a_i \l_i (\ov{a} -\langle a_i,\ov{a}\rangle a_i)$ (noting that this quantity belongs to the tangent space of $\partial\mathbb{S}^n_+$ at the point $a_i$) and summing for $i\in A_z$,  it holds that
\be\label{z*2} - c_6 \sum_{i\in A_z} \a_i^2 \frac{\n K_1(a_i) }{K_1(a_i)} (\ov{a} -\langle a_i,\ov{a}\rangle a_i)  - \frac{c_2}{2} \sum_{i, j\in A_z; j\neq i} \a_i \a_j \frac{\partial \e_{ij}}{\partial a_i}  (\ov{a} -\langle a_i,\ov{a}\rangle a_i) = o\big( \frac{1}{\mu_1} \big).\ee
Observe that, Proposition \ref{p:35} implies that $\a_i^2 = K(a_i)^{(2-n)/2} + O( \ln(\l_i) / \mu_1)$. Using Lemma \ref{derK} we get that
\begin{align}
 \sum_{i\in A_z} & \frac{\a_i^2}{K_1(a_i)}  \n K_1(a_i)   (\ov{a} -\langle a_i,\ov{a}\rangle a_i) \nonumber  \\
 &  =  \sum_{i\in A_z} \frac{1}{K_1(a_i)^{n/2}}  \n K_1(a_i)  (\ov{a} -\langle a_i,\ov{a}\rangle a_i) + O\Big( \frac{\ln\l_i}{\mu_1} d(a_i,\ov{a}) \Big) \nonumber \\
& =  -  \sum_{i\in A_z} \frac{1}{K_1(\ov{a})^{n/2}}  \n K_1(\ov{a})  ( a_i-\langle a_i,\ov{a}\rangle \ov{a}) + O( d(a_i,\ov{a})^2) + o\Big(\frac{1}{\mu_1}\Big) = o\Big(\frac{1}{\mu_1}\Big)\label{z*3} \end{align}
by using \eqref{rty1}. Finally, using Lemma \ref{lemderepsilon}, there holds
\be \frac{\partial \e_{ij}}{\partial a_i}  (\ov{a} -\langle a_i,\ov{a}\rangle a_i) + \frac{\partial \e_{ij}}{\partial a_j}  (\ov{a} -\langle a_j,\ov{a}\rangle a_j) \geq  c \e_{ij} . \label{z*4}\ee
Hence \eqref{z*2}-\eqref{z*4} contradict  Claim 5. Therefore the proof of Proposition \ref{comparable} is completed.
\end{pf}

\begin{rem}  We point out that for $\e = 0$,  each sequence of solutions $\o_k$ of $(\mathcal{P}_0)$ cannot be in $V(q,\ell,\tau)$ for $k $ large. Indeed the assumption that $ N \notin I$,  made in  the beginning of the proof,  is only used  in the proof of Claim 1.
\end{rem}

\subsubsection{ Location of blow up points and the speeds of concentration}
In the sequel we will consider the case where $\e >0$. In the next proposition we characterize  the location of blow up points and provide the   rate of the concentration parameters $\l_i$. Namely we prove

\begin{pro}\label{cp} Assume that $\e >0$.

$(a)$ Every interior concentration point $a_i$ converges to a critical point  $y_i$ of $K$ with  $\D K(y_i) < 0$,   $\l_i d( a_i, y_i ) $ is uniformly bounded and $y_i$ is an isolated simple blow up point. Moreover the concentration speed satisfies: $$ - c_4\frac{\D K(y_i) }{ \l_i ^2 K(y_i)} = 2 c_5  \e (1+o(1)) ,$$
where $c_4$ and $c_5$ are dimensional constants defined in Proposition \ref{p:33}.

$(b)$ Every boundary concentration point  $a_j$ converges to a critical point  $z_j$ of $K_1$ with  $\partial_{\nu} K  (z_j) > 0$. Furthermore, it holds :
$$  c_3 \frac{\partial_{\nu} K (z_j)}{K(z_j)}  \frac{1}{\l_j} = c_5 \, \e (1+o(1)),$$
where $c_3$ and $c_5$ are dimensional constants defined in Proposition \ref{p:33}.
\end{pro}

\begin{pf}  We start by proving Claim $(a)$. By Proposition \ref{comparable} we have that $I= I_b \cup I_{in}$. Hence, for each $i \in I_{in}$, it follows that,  $\mu_i:= \l_i^2$ and $\mu_1$ are of the same order. Hence using $(ii)$ of Lemma \ref{56} and $(E_i')$,  we derive that the set $D_1'$ has to be the empty set (indeed: if there exists $j\in D_1'$ then $a_j$ will converge to a critical point $y$ of $K$  and therefore $ | \D K(a_j) | > c >0$. Thus $(ii)$ of Lemma \ref{56} and $(E_j')$ are not compatible). Thus, from $(i)$ of Lemma \ref{56}, it follows that $a_i$ converges to a critical point $y_i$ of $K$  and we have $\l_i d( a_i , y_i ) \leq c$.  The sign of  $\D K(y_i)$ and the behavior of the $\l_i$ follow from $(E_i')$. Now assume that there exist  $i \neq j \in I_{in}$ such that $y_i = y_j$. Since $\l_k d( a_k , y_k ) \leq c$ for $k=i,j$ and $\l_i/\l_j$ is bounded from below and above, we derive that $\l_i\l_j d(a_i,a_j)^2$ is bounded and this is not compatible with  the fact that $\e_{ij}$ is small.  Hence the proof of the assertion $(a)$ is completed.\\
Concerning the second assertion, let $i\neq j \in I_b$, since $\l_i$ and $\l_j$ are of the same order we derive that $\l_k d(a_i,a_j) \to \infty$ as $\e \to 0$ for $k=i,j$  (by using the smallness of $\e_{ij}$). Using Lemma \ref{sum},  it holds
$$\frac{1}{\l_i} | \frac{\partial \e_{ij}}{\partial a_i} | \leq c \, \l_j d( a_i , a_j ) \e_{ij}^{n/(n-2)} \leq    \frac{c}{\l_i  d( a_i , a_j ) } \e_{ij} = o(\e_{ij}) = o(1/\mu_1).$$
Now using $(F_i')$ and the fact that $\l_i$ and $\mu_1$  are of the same order, we derive that $| \n K_1(a_i) | = o(1)$ for each $i\in I_b$ and therefore, $a_i$ has to converge to a critical point $z_i$ of $K_1$. To conclude the sign of $\partial K/ \partial \nu(z_i)$, we need the following claim \\
{\bf Claim A:}  $\e_{ij} = o(1/\mu_1)$ for each $i\neq j \in I_b$.\\
Arguing by contradiction, we assume that, there exists $i \in I_b$ such that $\sum_{j \neq i; j\in I_b} \e_{ij} \geq c / \l_i$. We remark that, in this case, the critical point $z_i$ (satisfying $a_i \to z_i$) has to satisfy $\# B_{z_i} \geq 2$. In fact, if not, i.e. $\# B_{z_i} =1$ we derive hat $| a_i - a_j | \geq c $ for each $j \neq i$ and therefore $\e_{ij}= O(1/(\l_i\l_j)^{(n-2)/2})$. Furthermore, we notice that Claim 5 of Proposition \ref{comparable} holds true (in fact, the proof relies on $(i)$ and this information is assumed in our case). Thus, arguing as in the end of the proof of Proposition \ref{comparable}, we derive a contradiction which implies the proof of our claim.\\
To achieve the proof of the proposition, we use Claim A and $(E_i')$ for $i\in I_b$.
\end{pf}

\subsubsection{Type of boundary blow up points }

In this subsection we study more carefully the boundary blow up points giving a precise characterization of \emph{isolated simple} blow up points and the \emph{non simple} blow up points.
\begin{lem}\label{nonsimple1}
Assume that all the $\mu_i$'s are of the same order and let $i_0 \neq  i_1 \in I_b$ be such that $d(a_{i_0}, a _{i_1}) := \min \{ d( a_i, a_j ) : i \neq j \in I_b \}$. Then we have  that:

\begin{enumerate}
  \item[(i)]
 Either  there exists a constant $c$ such that  $d( a_{i_0}, a _{i_1}) \geq c >0$ and in this case we have that  $\# B_z \leq 1$ for each critical point  $z$ of $K_1$,
  \item[(ii)] Or   $d(a_{i_0}, a _{i_1}) \to 0$ as $\e \to 0$ and in this case, there exists at least one critical point  $z$ of $K_1$ with $\# B_z \geq 2$ and  there   exist  positive constants $c, c', \tilde{c}$  such that  the following estimates hold:
   \begin{enumerate}
\item    $  \l_i ^{(n-2)/n} d(a_i, a_j ) \geq c  $ for each $i \neq j  \in I_b$.
\item $ \l_i ^{(n-2)/n} d(a_i, z ) \leq \tilde{c} $ for each $z$ such that $\# B_z \geq 2$ and for each $i \in B_{z}$.
\item   $ c \leq \l_i ^{(n-2)/n} d( a_i, a_j ) \leq c' $  for each $z$ such that $\# B_z \geq 2$ and for each $i \neq j  \in B_{z}$.
\end{enumerate}
\end{enumerate}
\end{lem}

\begin{pf}
The first assertion $(i)$ is immediate. Hence we will focus  on  the second assertion $(ii)$.\\
Since $\l_i$ and $\l_j$ are of same order for each $i,j \in I_b$, we derive that $\e_{ij} \leq c \e_{i_0i_1}$  for each $i,j \in I_b$.
Now, let $z$ be a critical point of $K_1$ such that $\# B_z \geq 2$ and let $k_0 \neq k_1 \in B_z$.

Observe that, for each $i\in B_z$ and  for each $j \in I_{in} \cup (I_b \setminus B_{z})$, it holds that $\e_{ij} \leq c/(\l_i \l_j)^{(n-2)/2} \leq c / \l_i^2$. Furthermore, for $j\in I_b$, it holds that $\l_j d( a_i , a_j ) \e_{ij}^{(n+1)/(n-2)} \leq c \e_{ij}^{n/(n-2)}$. Thus $(F_i)$ (defined in \eqref{Fib}) implies
\be\label{Fi''} \, \,   (F_i '' ) \quad  - \frac{c_2}{2} \sum_{j\in B_{z}; j\neq i} \a _j \frac{\partial\e _{ij}}{\partial a_i } -   c_6 \, \frac{ \a _i} {K(a_i)}{\n K_1}(a_i) = O\big( \frac{1}{\l_i} + \l_i \e_{i_0i_1}^{\frac{n}{n-2}} \ln (\e_{i_0i_1}^{-1}) \big).  \ee
Now, let $j_0$ be defined by: $| \n K_1(a_{j_0}) | :=\max\{ | \n K_1(a_{i}) |: i\in B_{z}\}.$

{\bf Claim 1:} There exists $M $ such that : $\displaystyle{ \frac{d( a_{j_0} , z )}{ \l_{j_0}} \leq \frac{M}{(\l_{i_0} d (a_{i_0} , a_{i_1} ) )^{n-1}}}$. \\
Arguing by contradiction, we assume that such a $M$ does not exist. Since the $\l_j$'s are of the same order and $z$ is nondegenerate, it follows that  : $ \e_{i_0i_1}^{(n-1)/(n-2)} = o( | \n K_1(a_{j_0}) | /\l_{j_0})$. Hence \eqref{Fi''} (with $i=j_0$) implies :
$$\frac{ | \n K_1(a_{j_0}) |}{\l_{j_0}} + \sum_{j\in B_{z}; j\neq j_0} O ( \e_{jj_0}^{(n-1)/(n-2)} ) = O\big( \frac{1}{\l_{j_0}^2} +  \e_{i_0i_1}^{n/(n-2)} \ln (\e_{i_0i_1}^{-1}) \big).$$
Since $ \e_{ij} \leq c \e_{i_0i_1}$  for each $i,j \in I_b$ and $ \e_{i_0i_1}^{(n-1)/(n-2)} = o( | \n K_1(a_{j_0}) | /\l_{j_0})$, then it follows that $\l_{j_0} | \n K_1(a_{j_0}) | $ is bounded. Therefore, we get that $\l_{j_0} d( a_{j_0} , z )$ is bounded. In addition, by the definition of $j_0$, we get that $ d( a_i , z ) \leq c d (a_{j_0} , z )$ for each $i \in B_{z}$. Thus, since the $\l_j$'s are of the same order, we derive that $\l_i d(a_i , z ) $ is bounded  for each $i\in B_{z}$. This implies that $\l_i \l_j d( a_i , a_j )^2$ is bounded which contradicts the smallness of $\e_{ij}$. Hence Claim 1 follows.

Note that, Claim 1 implies that
\be\label{xx13} \exists\, \,  M \mbox{ such that } \forall \, z \, (\mbox{with } \# B_z \geq 2) , \, \, \forall i \in B_z \, \, \mbox{ it holds } \, \,  \frac{d( a_{i} , z )}{ \l_{i}} \leq \frac{M}{(\l_{i_0} d( a_{i_0} , a_{i_1} ) )^{n-1}}.\ee

{\bf Claim 2:} There exists $c >0 $ such that : $\displaystyle{ \frac{d( a_{i_0} , a_{i_1} )}{ \l_{i_0}} \geq \frac{c}{(\l_{i_0} d( a_{i_0} , a_{i_1} ) )^{n-1}}}$.\\
Recall that $ d( a_{i_0} , a_{i_1} ) := \min\{ d( a_i , a_j ) : i\neq j \in I_b\}$. Let $z_0$ be such that $a_{i_0} \to z_0$ and let
$$ A'_{z_0} := \{ i \in B_{z_0} : \lim d( a_i , a_{i_0} ) / d( a_{i_0} , a_{i_1} ) = \infty\} \quad \& \quad A_{z_0} := B_{z_0} \setminus A'_{z_0}.$$
It follows that
\begin{itemize}
\item  $\forall \, \,  i \neq j \in A_{z_0}$, $d( a_i , a_j )$ and $d( a_{i_0} , a_{i_1} )$ are of the same order.
\item $\forall \, \,  i\in A_{z_0}$  and $\forall \, \,   j \in A'_{z_0}$, it holds : $ \lim d( a_i , a_{j} ) / d( a_{i_0} , a_{i_1} ) = \infty$.
\end{itemize}
Since the $\l_k$'s are of the same order, the second assertion implies that:
\be\label{mm1} \frac{1}{\l_i} | \frac{\partial \e_{ij}}{\partial a_i} | \leq   \frac{c}{\l_i}\frac{1}{(\l_i\l_j)^{(n-2)/2}} \frac{1}{d( a_i , a_j ) ^{n-1}} = o(\e_{i_0i_1}^{(n-1)/(n-2)}) \quad \forall i\in A_{z_0}\, , \forall \, j\in A'_{z_0}.\ee
Let  ${\bf \ov a}$ be defined as
\be\label{bara1}  {\bf \ov a} := \frac{b}{| b| } \quad \mbox{ where } \quad b \in \R^{n+1} \mbox{ satisfying }  \sum_{j\in A_{z_0}} (b - a_j) = 0.\ee
Note that \eqref{rty1} holds with $z=z_0$.
Multiplying $(F_i'')$ (defined in \eqref{Fi''}) by $\a_i( {\bf \ov a} - \langle a_i, {\bf \ov a}\rangle a_i)$  and summing over $i \in A_{z_0}$, we obtain (by using \eqref{mm1})
\begin{align}\label{mm2}
  - \frac{c_2}{2} \sum_{j\in A_{z_0}; j\neq i} \a_i \a _j \frac{\partial\e _{ij}}{\partial a_i } ( {\bf \ov a} - \langle a_i, {\bf \ov a}\rangle a_i) & -   c_6 \, \frac{ \a _i ^2} {K(a_i)}{\n K_1}(a_i)( {\bf \ov a} - \langle a_i, {\bf \ov a}\rangle a_i)\nonumber \\
  & = O\big( \frac{ d(  {\bf \ov a} , a_i ) }{\l_i} \big)+o\big( \l_i d(  {\bf \ov a} , a_i )  \e_{i_0i_1}^{(n-1)/(n-2)}  \big). \end{align}
Observe that Lemma \ref{lemderepsilon} gives the estimate of the first term. Now, using \eqref{alpha1} and Proposition \ref{cp}, we derive that
$$ \alpha_i^2 = \frac{1}{K(a_i)^{(n-2)/2}} + O(\e | \ln \e |).$$
Furthermore, using Claim $(i)$ of Lemma \ref{derK} (with $h={\bf \ov a}$) and \eqref{rty1},  \eqref{mm2} implies
\be\label{emna15} \sum_{i\neq j \in A_{z_0}} \e_{ij} \leq  c \sum_{i \in A_{z_0}} d( { \bf \ov a} , a_i )^2 +  c \sum_{i \in A_{z_0}} \frac{ d( {\bf \ov a} , a_i ) }{\l_i} + \sum_{i \in A_{z_0}} o\big( \l_i d(  {\bf \ov a} , a_i )  \e_{i_0i_1}^{\frac{n-1}{n-2}} + d(a_i,z) d({ \bf \ov a} , a_i ) \big) .\ee
Note that \eqref{xx13} implies that $ d(a_i,z) d(a_{i_0} , a_{i_1} ) \leq c \e_{i_0i_1}$. In addition, since $\l_i d( a_{i_0} , a_{i_1} )$ is very large  and $ d({ \bf \ov a} , a_i ) \leq c d( a_{i_0} , a_{i_1} )$ for each $i \in A_{z_0}$, we get that $ d({ \bf \ov a} , a_i ) / \l_i = o( d( a_{i_0} , a_{i_1} )^2)$. Furthermore, using the fact that $\l_i d({\bf \ov a} ,  a_i ) \e_{i_0i_1}^{1/(n-2)} \leq c $ for each $i \in A_z$, then \eqref{emna15} implies that $ \e_{i_0i_1} \leq c\,  d( a_{i_0} , a_{i_1} )^2$ which implies Claim 2.
\medskip

Now, we are in position to prove Assertion $(ii)$ of the lemma. In fact:\\
 Claim 2 implies that $ \l_{i_0}^{(n-2)/n} d( a_{i_0} , a_{i_1} ) \geq c$ and therefore, since all the $\l_j$'s are of the same order and $d( a_{i_0} , a_{i_1} )  = \min \{ d( a_{i} , a_{j} ) : i \neq j \in I_b\}$, we get that \be\label{xx12} \l_i^{(n-2)/n} d( a_{i} , a_{j} )  \geq c \quad \mbox{ for each } i\neq j \in I_b.\ee
Hence the proof of Assertion $(a)$ (in particular, the first inequality of the assertion $(c)$) follows.\\
Furthermore, let $z$ be such that $\# B_z \geq 2$, and let $i \in B_z$, using \eqref{xx13}, we get
\be\label{xx11} \l_i^{\frac{n-2}{n}} d( a_{i} , z ) =  \l_i^{\frac{2n-2}{n}} \frac{d( a_{i} , z )}{ \l_i} \leq \frac{M \l_i^{(2n-2)/n}}{(\l_{i_0} d( a_{i_0} , a_{i_1} ) )^{n-1}} \leq \frac{C }{(\l_{i_0}^{(n-2)/n}  d( a_{i_0} , a_{i_1} ) )^{n-1}} \leq c \ee (by using \eqref{xx12}) which gives the proof of Assertion $(b)$.

Now, for $z$ such that $\# B_z \geq 2$ and  for $i\neq j \in B_{z}$,  using the fact that $| a_{i} - a_{j} | \leq | a_{i} - z | + |  a_{j} -z |$ and applying \eqref{xx11}, we derive the second inequality of Assertion $(c)$ which completes the proof of Assertion $(c)$.\\
Note that, once Assertion $(c)$ is proved, we conclude that $A'_{z} = \emptyset$ for each $z$ such that $\# B_z \geq 2$.
The proof of the lemma is thereby completed.
\end{pf}

Note that from \eqref{xx12}, it follows that  (we can be more precise if $d( a_i , a_j ) \geq c $)
\be\label{eps5}
\e_{ij} \leq \frac{c}{\l_i ^{2(n-2)/n}} \, \, \mbox{ for each } i\neq j \in I_b.\ee
Furthermore, Assertion $(c)$ of Lemma \ref{nonsimple1} implies that, for each $z$ such that $\# B_z \geq 2$, all the $d( a_i , a_j )$'s, for $i\neq j \in B_z$, are of the same order. Hence $B_z = A_z$ for each $z$.

  \begin{lem}\label{bz1}
  Let $z$ be such that $ B_z =\{i\}$. Then it holds that $ \l_i d( a_i , z ) \leq c$.
  \end{lem}

\begin{pf}
We need to use Equation $(F_i)$ (introduced in \eqref{Fib}) but in this equation, it appears one term in $R_{a_i}$ which we cannot control. This term is $\sum _{k\neq j}\e_{kj}^{n/(n-2)} \ln(\e_{kj}^{-1})$ (by using \eqref{eps5}, it is enough to get $\gamma > n/(n-2)$ instead of $n/(n-2)$). For this reason, we will repeat the proof of the equation $(F_i)$  in our special case, that is $| a_j - a_i | \geq c $ for each $j \neq i$. In fact, we will follow the proof of the third assertion of Proposition 3.1 of \cite{AB20b} (by taking $\o = 0$). Multiplying the equation satisfied  by $u_\e$ by $(1/\l_i) (\partial \d_i/\partial a_i^k)$ we get
\be\label{nbv2}  \langle u_\e ,  \frac{1}{\l_i} \frac{ \partial \d_i}{\partial a_i^k}\rangle = \int_{\mathbb{S}^n_+}  K {u}_\e^{p-\e}  \frac{1}{\l_i} \frac{ \partial \d_i}{\partial a_i^k} \quad \mbox{ with } \quad p:= \frac{n+2}{n-2}.\ee
Since $ v_\e $ satisfies \eqref{eq:V0} and  $a_i \in \partial \mathbb{S}^n_+$, the first term is $\sum_{j \neq i} O(\e_{ij})$ which is $o(1/\l_i ^2)$ (by using the fact that $d(a_i,a_j) > c >0$ for each $j\neq i$). Furthermore, the other term is estimated as
\be\label{nbv1} \int_{\mathbb{S}^n_+}  K {u}_\e^{p-\e}  \frac{1}{\l_i} \frac{ \partial \d_i}{\partial a_i^k}  = \int_{\mathbb{S}^n_+}  K \ov{u}_\e^{p-\e}  \frac{1}{\l_i} \frac{ \partial \d_i}{\partial a_i^k}  + (p - \e ) \int_{\mathbb{S}^n_+}  K \ov{u}_\e^{p-1-\e}  \frac{1}{\l_i} \frac{ \partial \d_i}{\partial a_i^k} v_\e + O(\| v_\e \|^2) \ee
where $\ov{ u } := u_\e - v_\e$. Furthermore, using the behavior of $\l_i$ (given in Claim $(b)$ of Proposition \ref{cp}) and \eqref{eps5}, we derive that $R^2(\e,a,\l)$ (which is introduced in Proposition \ref{eovv}) is $O(1/\l_i^2)$ which implies that $\| v_\e \|^2 =  O(1/\l_i^2)$. In addition, observe that
\begin{align*} \Big|  \ov{u}^{\frac{4}{n-2} - \e } \frac{1}{\l_i} \frac{ \partial \d_i}{\partial a_i^k} - (\a_i \d_i )^{\frac{4}{n-2} - \e } \frac{1}{\l_i} \frac{ \partial \d_i}{\partial a_i^k} \Big| & \leq c | \ov{u} ^{\frac{4}{n-2} - \e }  - (\a_i \d_i )^{\frac{4}{n-2} - \e } | \d_i \\
& \leq c \sum_{j \neq i}  \begin{cases} (\d_i \d_j)^{\frac{n+2}{2(n-2)}} \mbox{if } n \geq 6,\\
\d_i ^\frac{4}{3} \d_j  + \d_i \d_j ^\frac{4}{3}  \mbox{if } n = 5. \end{cases} \end{align*}
Hence, following the proof of Eq (23)  of \cite{AB20b}, it follows that the second term of \eqref{nbv1} is $O(R(\e,a,\l) \| v_\e \| )$ and therefore it is also $O(1/\l_i^2)$.
 For the first integral of the right hand side of \eqref{nbv1}, using the behavior of $\l_i$ (given in Proposition \ref{cp}),  it holds that $$\int K \ov{u}_\e^{p-\e}  \frac{1}{\l_i} \frac{ \partial \d_i}{\partial a_i^k} = \a_i^{p-\e} \int K \d_i^{p-\e}  \frac{1}{\l_i} \frac{ \partial \d_i}{\partial a_i^k}  + \sum_{j\ne i}O(\e_{ij}) = \a_i^{p-\e}\frac{c}{\l_i} \frac{\partial K}{\partial x^k}(a_i) + O(\frac{1}{\l_i^2}). $$ Thus,  \eqref{nbv2} becomes $$  \frac{ \n K_1 (a_i) } {\l_i } = O(\frac{1}{\l_i^2} )$$ which implies the result (since the critical points of $K_1$ are  non-degenerate).
\end{pf}

\subsection{Proof of  Theorem \ref{th:t1}}

The assertions $(a)$ and $(b)$ follow from Proposition \ref{cp} by choosing
\be \label{kappa12} \kappa_1(n):= c_4/ (2c_5) \quad \& \quad \kappa_2(n):= c_3/ c_5 \ee where the constants $c_3$, $c_4$ and $c_5$ are defined in Proposition \ref{p:33}.

 It remains to see the last one.
Let  $z$ be a critical point of $K_1$ such that $\# B_z \geq 2$ and let $i_0$ and $i_1$ be such that $d( a_{i_0} , a_{i_1} ) := \min \{ d( a_i , a_j ): i\neq j \in B_z\}$.   Let $b_i$ be defined by \eqref{bi}. From Assertion $(c)$ of Lemma \ref{nonsimple1} and the behavior of the $\l_i$'s (given in Claim $(b)$ of the theorem), we derive that $0 < c \leq | b_{i} - b_{j} |  \leq c' $ for each $i \neq j$. In addition, using the behavior of $\l_i$ and using Assertion $(b)$ of Lemma \ref{nonsimple1},  we derive that  $ | b_i | \leq C$ for each $i \in B_z$. Hence, it follows that
\be \label{mm9} | b_i | \leq C \quad \mbox{ and } \quad | b_i - b_j | \geq c   \qquad \forall \, \, i\neq j \in B_z.\ee
Furthermore  observe that : $2(1-\cos(d(a_i,a_j))) = | a_i - a_j |^2$  (seen as two points of $\R^{n+1}$) and we notice that, if $d(a,b)$ is small, it holds that $d(a,b) = | a-b| (1+o(1))$. Therefore, from the definition of $\e_{ij}$ (see \eqref{eppe}) and the fact that $\l_i$ and $\l_j$ are of he same order, it follows that
\be\label{mm10}  \frac{\partial \e_{ij}}{\partial a_i} = \frac{n-2}{4}\l_i \l_j (a_j-a_i) \e_{ij}^{n/(n-2)} = \frac{2^{n-2} (n-2)}{(\l_i \l_j)^{(n-2)/2}} \frac{a_j - a_i}{ | a_j - a_i |^n}  +O\Big(\l_i \e_{ij}^{\frac{n+1}{n-2}}\Big).\ee
In addition, using \eqref{alpha1}, it holds that
\begin{align}\frac{\a_j}{\a_i} & = \Big(\frac{K_1(a_i)}{K_1(a_j)}\Big)^{(n-2)/4} + O(\e | \ln \e | ) \nonumber \\
&  = 1 + O(\e | \ln \e | + d( a_i , z )^2 + d( a_j ,z )^2 ) = 1 + O(\e | \ln \e | ). \label{mm11}  \end{align}

Hence, for each $e\in \partial \mathbb{S}^n_+$ such that $ c \leq d( e , z ) \leq c' $, using Lemma \ref{derK}  and \eqref{mm10}, \eqref{mm11}, $(F_i'')$ implies
\begin{align}  -\frac{c_2}{2} \sum_{j\in B_z; j\neq i} & \frac{2^{n-2} (n-2)}{(\l_i \l_j)^{(n-2)/2}} \frac{1}{ | a_j - a_i |^n} \langle a_j - a_i, e - \langle e,a_i\rangle a_i\rangle \nonumber\\
&  -  c_6  \frac{1}{  K_1(z)} D^2 K_1(z) \big(a_i - \langle a_i, z \rangle z, e - \langle e, z \rangle z\big) + O(  d ( a_i , z )^2 ) = o_{\e}(\e ^{(n-2)/n}).\label{mm13}\end{align}

Note that it is easy to see that: $| a_i - z |^2 = 2 (1- \langle a_i , z \rangle ) $ and therefore, it holds that  $ \langle a_i , z \rangle = 1 + O(\e^{2(n-2)/n})$ (by using $(b)$ of Lemma \ref{nonsimple1} and the fact that $d(a_i,z) = | a_i-z| (1+o(1))$) which implies that $b_i$ (defined by \eqref{bi}) satisfies $b_i= \e^{(2-n)/n} (a_i - z) + O(\e^{(n-2)/n})$. Thus, using the behavior of the $\l_i$'s (see Claim $(b)$) and the change of variables \eqref{bi}, and multiplying \eqref{mm13} by $\e^{(2-n)/n}$, it follows that
$$ - \sum_{j\in B_z; j\neq i} \frac{n-2}{ | \ b_j - b_i | ^n} \langle b_j - b_i ,  e - \langle e, z \rangle z \rangle -   D^2 K_1(z) \big(b_i , e - \langle e, z \rangle z\big) = o_\e(1), $$
by choosing $\kappa _3(n)$ equal to
\be \label{kappa3} \kappa_3(n) := 2^{(n-3)/n} (c_2/ c_6)^{1/n} (c_5/c_3)^{(n-2)/n}. \ee

Finally, using \eqref{mm9}, we derive that , for each $i \in B_z$, $b_i$ converges to $\ov{b}_i$ satisfying
\begin{align*}
&  \ov{b}_i  \in T_z(\partial \mathbb{S}^n_+) \, \, , \, \,  | \ov{b}_i | \leq c \, \, , \, \, | \ov{b}_j - \ov{b}_i | \geq c \, \, \forall \, \, i\neq j  \quad \mbox{ and }\\
&  \sum_{j\in B_z; j\neq i} \frac{n-2}{ | \ov{b}_j - \ov{b}_i | ^n} \langle \ov{b}_j - \ov{b}_i ,  e - \langle e, z \rangle z \rangle + D^2 K_1(z) \big(\ov{b}_i, e - \langle e, z \rangle z\big) = 0 \end{align*}
which means that $(\ov{b}_1,\cdots,\ov{b}_m)$ is a critical point of $\mathcal{F}_{z,m}$ (where $m:= \# B_z$).\\
Thus the proof of the theorem is completed.




\section{The non zero weak limit case}

\begin{pro}\label{lambdaepsilonw}
Let $ u_\e$  be an energy bounded solution of $(\mathcal{P}_\e)$ which blows up. We assume that  there exists a positive solution $\o$ of $(\mathcal{P}_0)$ such that $u_\e \rightharpoonup \o$ (but $u_\e \nrightarrow \o$). Then there exist $q$ and $\ell$ such that $u_\e$ can be written as
 \begin{equation}\label{kkk111} u_\e:= u_{\a,\b} + \sum_{ i = 1}^q \a_i \d_i + \sum_{i  =q+1}^{q+\ell} \a _j \varphi _j + v_\e \in V(\o, q,\ell ,\tau )\quad \mbox{ with } v_\e\in E_{\o,a,\l}^\perp. \end{equation}
\end{pro}

In this section, we will adapt the program done in Section \ref{ss3} and we will present the contribution of $\o$ in the expansions proved in the previous section. In the sequel, we will collect some properties which are satisfied by the parameters and the solution $u_\e$. Our first result concerns the function $v_\e$ to show that it does not have any contribution in the phenomenon.

\begin{pro}\label{eovvw}
Let the  remainder term  $v_\e$  be defined as in Proposition \ref{lambdaepsilonw}. Then there holds:
$$ \| v_\e \| \leq c  R(\e,a,\l) + \sum \xi(\l_i)      \qquad  \mbox{ where }\quad  \xi(\l_i) :=  \begin{cases} 1/\l_i ^{(n-2)/2} \mbox{ if } n \leq 5\\
\ln(\l_i)^{2/3} / \l_i^2 \mbox{ if } n = 6\end{cases} $$ and where
$ R(\e,a,\l)$ is defined in Proposition \ref{eovv}.
\end{pro}
\begin{pf}
Recall that $\o$ is a solution of $(\mathcal{P}_0)$ (not necessary non degenerate), thus we decompose $H^1(\mathbb{S}^n_+)$ as
\be \label{vvv*1}  H^1(\mathbb{S}^n_+) := N_-(\o) \oplus H_0(\o) \oplus N_+(\o) \, \, \quad \mbox{with} \, \, \,  H_0(\o):= span\{\o\} \oplus N_0(\o)\ee
 where  $N_-(\o)$, $N_0(\o)$ and $ N_+(\o)$ are respectively the space of negativity, of nullity and of positivity  of the quadratic form $Q_\o$ (defined by \eqref{qw}) in $span\{\o\} ^\perp$. Note that these spaces are orthogonal spaces with respect to $\langle.,.\rangle$ and the associated bilinear form $B_\o(.,.)$ ($:= \int_{\mathbb{S}^n_+}  K \o^{p-1} . . $). Furthermore, the sequence of the eigenvalues (denoted by $(\s_i)$) corresponding to $Q_\o$ satisfies $ \s_i \nearrow 1$. Therefore, there exists a constant $c >0$ such that
\begin{equation}\label{qw+-}
Q_\o(h) \leq -c \, \| h \| ^2 \quad \forall \, \,  h \in N_-(\o) \quad ; \quad Q_\o(h) \geq c \, \| h \| ^2 \quad \forall \, \, h \in N_+(\o).\end{equation}
Using \eqref{vvv*1},  $v_\e $ can be decomposed as  follows
\be\label{emna13} v_\e := v_\e^- + v_\e^0 + v_\e^+ \quad \mbox{ where } v_\e ^0 \in  H_0(\o) \, \, ; \, \, v_\e ^- \in  N_-(\o) \mbox{ and }  v_\e ^+ \in  N_+(\o) .\ee
Since $v_\e \in E_{w,a,\l}^\perp$, using Lemma \ref{Qwalpositive}, we get that $ v_\e^0 = o( \| v_\e \| )$. For the other parts, note that $v_\e^+$ and $v_\e^-$ are not necessarily in $E_{a,\l}^\perp$ but they are in $H_0(\o)^\perp$. \\
Now we will focus on estimating $v_\e^-$. Taking $u=u_\e$ and $h=v_\e^-$ in \eqref{www1} and using \eqref{ti6}, we derive that
\begin{align}
 \sum \a_j \langle \varphi_j , v_\e^- \rangle +  \langle u_{\a,\b} , v_\e^- \rangle + \| v_\e^- \|^2 & = \int  K u_\e ^{p-\e} v_\e^- \nonumber\\
& =   \int K \ov{u}_\e ^{p-\e} v_\e^- + p \int K \ov{u} _\e ^{p-1-\e} v v_\e^- + o(\| v_\e\| \| v_\e^-|) \label{emna1}\end{align} where $\ov{u}_\e := u_\e - v_\e.$
Using \eqref{ti5}, \eqref{alphai0} and the fact that $v_\e^-$ is in a finite dimensional space  (which implies that $\| . \|_\infty \leq c \| . \|$), we derive that
\begin{align}
 \int K \ov{u} _\e ^{p-1-\e} &  v_\e v_\e^-  = \sum_{i=1}^N \int K (\a_i \varphi_ i) ^{p-1-\e} v_\e  v_\e ^- + \int K u_{\a,\b}^{p-1-\e} v_\e  v_\e ^- + o(\| v_\e\| \| v_\e ^-\|)\nonumber \\
 & = \sum_{i=1}^N O\big( | v_\e ^- |_\infty \int  \d_ i ^{p-1} | v_\e | \Big)  +  \int K \o^{p-1} ( v_\e ^0 +v_\e ^- + v_\e^+ ) v_\e^- + o(\| v_\e\| \| v_\e ^-\|) \nonumber \\
& =  \int K \o^{p-1} ( v_\e ^-  )^2 + o(\| v_\e\| \| v_\e ^-\|)  \label{emna2}
 \end{align} where we have used the orthogonality of $v_\e ^-$, $v_\e^0$ and $ v_\e ^+$ with respect to $\int K\o^{p-1} ..$ . Thus \eqref{emna1} and  \eqref{emna2} imply that
\be\label{emna3}  - Q_\o (v_\e^-) + o(\| v_\e\| \| v_\e ^-\|) =  \sum \a_j \langle \varphi_j , v_\e^- \rangle +  \langle u_{\a,\b} , v_\e^- \rangle  - \int K \ov{u}_\e ^{p-\e} v_\e^- := \ell (v_\e^-) . \ee
Observe that, using $ \| v_\e^- \| _\infty \leq c \| v_\e ^- \|$, it holds that
\be\label{emna11}
| \langle \varphi_j , v_\e^- \rangle | \leq c  \int \d_i ^p | v_\e^- |  \leq c   \| v_\e^- \|_\infty  \int \d_i ^p  \leq c  \frac{ \| v_\e^- \| }{\l_i ^{(n-2)/2}} \quad ; \quad  \int \d_i  | v_\e^- |  \leq c\frac{ \| v_\e^- \| }{\l_i ^{(n-2)/2}} . \ee
\begin{align} \int K \ov{u}_\e ^{p-\e} v_\e^- & = \int K {u}_{\a,\b} ^{p-\e} v_\e^- + \sum O\Big( \int u_{\a,\b} ^{p-1} \d_i | v_\e^-| + \int \d_i ^{p-\e} | v_\e^-| \Big) \nonumber \\
& = \int K {u}_{\a,\b} ^{p} v_\e^-  +  O\Big( \| v_\e^- \| \Big( \e + \sum \frac{1}{\l_i ^{(n-2)/2}}\Big)\Big) . \label{emna4}\end{align}
Therefore, combining \eqref{qw+-}, \eqref{emna3}-\eqref{emna4} with \eqref{ualphbeta2}, we get
\be\label{emna5}  c \, \| v_\e ^- \| ^2 \leq - Q_\o (v_\e ^- ) \leq c \, \| v_\e ^-  \| \Big( \e + \sum \frac{1}{\l_i ^{(n-2)/2}} + o( \| v_\e \| ) \Big) .\ee
It remains to estimate the $v_\e ^+$-part. Recall that $ v_\e ^+ \in N_+(\o)$ and it is not necessarily in $E_{a,\l}^\perp$. Note that Eq \eqref{emna1} holds with $v_\e ^+$ instead of $v_\e ^-$. Observe that, using \eqref{ti5} and \eqref{alphai0}, it holds
\begin{align}
 \int K \ov{u} _\e ^{p-1-\e} v_\e v_\e^+ & = \sum_{i=1}^N \int K (\a_i \varphi_ i) ^{p-1-\e} v_\e  v_\e ^+ + \int K u_{\a,\b}^{p-1-\e} v_\e  v_\e ^- + o(\| v_\e\| \| v_\e ^+\|)\nonumber \\
 & = \sum_{i=1}^N   \int  \d_ i ^{p-1}  v_\e   v_\e ^+  +  \int K \o^{p-1} ( v_\e ^0 +v_\e ^- + v_\e^+ ) v_\e^+ + o(\| v_\e\| \| v_\e ^+\|) \nonumber \\
& =   \sum_{i=1}^N   \int  \d_ i ^{p-1}   (v_\e ^+)^2  + \int K \o^{p-1} ( v_\e ^+  )^2 + o(\| v_\e\| \| v_\e ^+\|) , \label{emna22}
 \end{align}
 by using the fact that $ \| v_\e^- \| _\infty \leq c \| v_\e ^- \|$ and $ \| v_\e^0 \| _\infty \leq c \| v_\e ^0 \|$ and the orthogonality of $v_\e ^-$ (respectivement $v_\e^0$) and $ v_\e ^+$ with respect to $\int K\o^{p-1} ..$. Hence we obtain
 \be\label{emna9} Q_{\o,a,\l} (v_\e^+) +  o(\| v_\e\| \| v_\e ^+\|) =  - \sum \a_j \langle \varphi_j , v_\e^+ \rangle -  \langle u_{\a,\b} , v_\e^+ \rangle  +  \int K \ov{u} _\e ^{p-\e} v_\e^+ := - \ell (v_\e^+).\ee
Using Lemma \ref{Qwalpositive}, we derive that $Q_{\o,a,\l} (v_\e^+) \geq c \| v_\e^+\|^2 + o(\| v_\e\|^2)$, hence it remains to estimate the linear part $\ell(v_\e^+)$. In fact, using  \eqref{ti5}, we have
\begin{align} \int K \ov{u}_\e ^{p-\e} v_\e^+ = & \sum \int_{\mathbb{S}^n_+} K ( \a_i \varphi_i)^{p-\e} v_\e^+ + \int_{\mathbb{S}^n_+} K u_{\a,\b}^{p-\e} v_\e ^+ \nonumber\\
&  + O\Big( \begin{cases} \sum_{i\neq j} \int  (\d_i \d_j)^{p/2} | v_\e ^+|  + \sum \int (\d_i \o )^{p/2}  | v_\e^+ |  & \mbox{ if } n \geq 6 \\
 \sum_{i\neq j} \int  \d_i^{4/3} \d_j  | v_\e ^+|  + \sum \int ( \d_i ^{4/3} \o + \d_i \o^{4/3} ) | v_\e^+ |  & \mbox{ if } n = 5 \end{cases} \Big) . \label{ff1} \end{align}
Concerning the remainder terms, using Lemma \ref{somest}, it follows that 
\begin{align}
&  \int  (\d_i \d_j)^{p/2} | v_\e ^+| \leq \| v_\e^+\| \Big( \int (\d_i \d_j )^{n/(n-2)} \Big)^{(n+2)/(2n)} \leq c  \| v_\e^+\|  \e_{ij}^{\frac{n+2}{2(n-2)}} \ln ^{\frac{n+2}{2n}} (\e_{ij}^{-1}) , \\
& \int_{\mathbb{S}^n_+} \d_i ^{p/2}  | v_\e^+ |  \leq \| v_\e^+\| \Big( \int_{\mathbb{S}_n^+} \d_i ^{n/(n-2)} \Big)^{(n+2)/(2n)} \leq c  \| v_\e^+\|  \frac{\ln ^{(n+2)/(2n)} (\l_i)}{\l_i ^{ (n+2)/4}} ,\\
&  \int_{\mathbb{S}^5_+}  \d_i^\frac{4}{3} \d_j  | v_\e ^+| =   \int  \d_i^\frac{1}{3} (\d_i \d_j)  | v_\e ^+|  \leq c \| v_\e^+\| \Big( \int (\d_i \d_j )^\frac{5}{3} \Big)^\frac{3}{5} \leq c  \| v_\e^+\|  \e_{ij} \ln ^\frac{3}{5} (\e_{ij}^{-1}) ,\\
& \int_{\mathbb{S}^5_+}  \d_i ^{4/3}  | v_\e^+ | \leq \| v_\e^+\| \Big( \int \d_i ^{40/21} \Big)^{7/10} \leq c  \| v_\e^+\| \frac{1}{\l_i ^{3/2}} ,\\
& \int_{\mathbb{S}5_+}  \d_i  | v_\e^+ |  \leq \| v_\e^+\| \Big( \int \d_i ^{10/7} \Big)^{7/10} \leq c  \| v_\e^+\| \frac{1}{\l_i ^{3/2}}.
\end{align}
For the other integrals, we recall that  $ v_\e ^+ \in N_+(\o)$ and therefore it follows that $ v_\e ^+ \in H_0(\o)^\perp$. Hence using \eqref{ualphbeta2} we derive that
\be\label{emna7} -  \langle u_{\a,\b} , v_\e^+ \rangle +  \int_{\mathbb{S}^n_+} K u_{\a,\b}^{p-\e} v_\e ^+ = -  \langle u_{\a,\b} , v_\e^+ \rangle +  \int_{\mathbb{S}^n_+} K u_{\a,\b}^{p} v_\e ^+  + O(\e \, \| v_\e^+\| ) = O(\e \, \| v_\e^+\| ) .\ee
Furthermore, using Lemma \ref{lowerL2}, it holds that
\begin{align} \int_{\mathbb{S}^n_+} K ( \a_i \d_i)^{p-\e} v_\e^+ = & \frac{\a_i^{p} K(a_i)}{\l_i ^{\e(n-2)/2}}  \int_{\mathbb{S}^n_+} \d_i ^{p} v_\e^+ + O\Big( \int_{\mathbb{S}^n_+} | K(x) - K(a_i) |   \d_i ^{p} | v_\e^+ | \nonumber \\
& + \e \int \d_i ^p \ln [2+(\l_i ^2 -1)(1-\cos d(a_i,x))] | v_\e^+ |  + \e \| v_\e^+\|    \Big) \nonumber\\
& = \frac{\a_i^{p} K(a_i)}{\l_i ^{\e(n-2)/2}}  \langle \d_i , v_\e^+ \rangle  + O\Big(   \| v_\e^+ \| \Big( \e + \frac{| \n K(a_i) | }{\l_i} + \frac{1}{\l_i ^2}\Big) \Big).\label{emna12} \end{align}
Thus, for $i \leq q$, (that is $a_i \in \partial \mathbb{S}^n_+$), we have $\varphi _i = \d_i$ and therefore we get
\begin{align} \label{emna8} -  \a_i \langle \d_i , v_\e^+ \rangle  & + \int_{\mathbb{S}^n_+} K ( \a_i \d_i)^{p-\e} v_\e^+ \\
& = \a_i \Big(   \frac{\a_i^{p} K(a_i)}{\l_i ^{\e(n-2)/2}} - 1 \Big) \langle \d_i , v_\e^+ \rangle  + O\Big(   \| v_\e^+ \| \Big( \e + \frac{| \n K(a_i) | }{\l_i} + \frac{1}{\l_i ^2}\Big) \Big) .\nonumber \end{align}
However, for $i \geq q+1$ (that is $a_i \in \mathbb{S}^n_+$), we have
$$ \int_{\mathbb{S}^n_+} K \varphi_i^{p-\e} v_\e^+ = \int_{\mathbb{S}^n_+} K \d_i^{p-\e} v_\e^+ + O\Big( \int_{\mathbb{S}^n_+}  \d_i^{p-1} | \varphi _i - \d_i | | v_\e^+ | \Big).$$
The first term is estimated in \eqref{emna12} and the second one satisfies
$$ \int_{\mathbb{S}^n_+}  \d_i^{p-1} | \varphi _i - \d_i | | v_\e^+ | \leq c \| v_\e ^+ \| \Big(  \int_{\mathbb{S}^n_+}  \d_i^{8n/(n^2-4)} | \varphi _i - \d_i |^{2n/(n+2)} \Big) ^{(n+2)/(2n)}$$
The estimate of the last integral depends on the dimension $n$. In fact,  for $n \geq 6$, we have $2n/(n+2) \geq n/(n-2)$ and  using the fact that $| \varphi _i - \d_i | \leq c \, \min (\d_i ; 1/(\l_i d_i ^2)^{(n-2)/2})$ (see Lemma \ref{lem:varphi}), we derive that
$$  \int_{\mathbb{S}^n_+}  \d_i^{\frac{8n}{n^2-4}} | \varphi _i - \d_i |^{\frac{2n}{n+2}}  \leq c | \varphi _i - \d_i |_\infty ^{n/(n-2)}  \int_{B(a_i,d_i)}  \d_i^{\frac{n}{n-2}} +  \int_{\mathbb{S}^n_+\setminus B(a_i,d_i)}  \d_i^{\frac{2n}{n-2}}  \leq c \frac{\ln(\l _i d_i ) }{ (\l_i d_i )^n }
$$
and for $n \leq 5$, we have $ 8n/(n^2-4) > n/(n-2)$ and therefore it holds
$$  \int_{\mathbb{S}^n_+}  \d_i^{8n/(n^2-4)} | \varphi _i - \d_i |^{2n/(n+2)}  \leq c | \varphi _i - \d_i |_\infty ^{2n/(n+2)}  \int_{\R^n_+}  \d_i^{8n/(n^2-4)}
  \leq c \frac{ 1 }{ (\l_i d_i )^{2n(n-2)/(n+2)} }.$$
Finally, using the fact that $v_\e \in E_{w,a,\l}^\perp$, as in the computations of \eqref{emna11}, we derive that
$$ \langle \varphi_i , v_\e^+ \rangle = \langle \varphi_i , v_\e \rangle - \langle \varphi_i , v_\e^0 \rangle - \langle \varphi_i , v_\e^- \rangle = O\Big( \| v_\e \| / \l_i ^{(n-2)/2} \Big) \quad \forall \, \, 1 \leq i \leq q+\ell.$$
Thus, using \eqref{qw+-},  \eqref{emna9} becomes
\be\label{emna133} c \| v_\e^+\| ^2 \leq Q_{\o,a,\l}(v_\e^+) \leq o( \| v_\e\| ^2) + O\Big( \| v_\e \| \Big( R(\e,a,\l) +  \sum \xi(\l_i) \Big) \Big).\ee
Thus using \eqref{emna13},  \eqref{emna5} and   \eqref{emna133} we derive that
$$ \| v_\e\| ^2 = \| v_\e^+\|^2 + \| v_\e^-\|^2 + o( \| v_\e\|^2) \leq c  \| v_\e \| \Big(R(\e,a,\l) +  \sum \xi(\l_i) \Big)$$
which completes the proof.
\end{pf}


\subsection{Balancing conditions for  blow ups with residual mass }

In this subsection we prove various balancing conditions for the parameters of concentration in the case where the blowing up solution $u_{\e}$  of $(\mathcal{P}_{\e})$ decomposes into a solution of the problem $(\mathcal{P})$ plus a sum of  interior and or boundary bubbles. The following propositions are quoted from \cite{AB20b} (Propositions 3.1 and 3.2) by taking $u_{\a,\b}$ instead of $\o$. Note that in our case we have $\n I_\e(u_\e) =0$. \\
We start by estimating the gluing parameters $\a_i$. Namely we have:
\begin{pro}\label{p:35w} \cite{AB20b}
For each $1\leq i \leq N$, it holds:
$$ | 1- \l_i^{-\e\frac{n-2}{2}} \a _i^{\frac{4}{n-2}}K(a_i) | = O\Big( R_{\a_i} + \frac{ 1 }{\l_i^{(n-2)/2}} + \sum \frac{ \ln \l_k }{\l_k^{n/2}} \Big)$$
where $ R_{\a_i} $ is defined in Proposition \ref{p:35}.
\end{pro}
Next we provide  balancing conditions involving the rate of the concentration $\l_i$ and the  mutual interaction of bubbles $\e_{ij}$.
\begin{pro}\label{p:33w}
  For $\e$ small enough, the following equations  hold:
  \begin{align*}
   - \frac{c_2} {2}\sum _{j\neq i; j\leq q} \a _j \l _i\frac{\partial \e _{ij}}{\partial \l _i}   -  \frac{\a_i }{K(a_i)} & \Big[ \frac{c_3}{\l _i}\frac{\partial K}{\partial \nu}(a_i) -  c_5 K(a_i) \e\Big]  + c_2\frac{n-2}{4} (1+o(1)) \frac{ \o(a_i)}{\l_i ^{(n-2)/2}}\\
   & = O (\frac{1}{\l_i ^2} + \sum_{j > q} \e_{ij} + R_{1,\l} +  \sum \frac{\ln \l_k }{\l_k ^{n/2}}   ) \, \qquad  (\mbox{for } i \leq q) \end{align*}
\begin{align*}   - {c_2} \sum _{j\neq i} \a _j \l _i\frac{\partial \e _{ij}}{\partial \l _i}  & + c_2  \frac{n-2}{2}  \sum_{j=q+1}^p\a_j  \frac{H(a_i,a_j)}{(\l_i\l_j)^{(n-2)/2}}  +  {\a_i} \Big(c_4 \frac{\Delta K(a_i)}{\l_i^2 K(a_i)} + 2c_5 \e\Big)\\
 &  + c_2\frac{n-2}{2} (1+o(1)) \frac{ \o(a_i)}{\l_i ^{(n-2)/2}} = O (  R_{2,\l} + \sum \frac{\ln \l_k }{\l_k ^{n/2}}  ) \, \qquad  (\mbox{for } i > q) \end{align*}
 where $ R_{1,\l} $,  $ R_{2,\l} $, $ c_3 $, $ c_2$ and $ c_5 $ are defined in Proposition \ref{p:33}.
\end{pro}

Finally we provide the following balancing conditions involving  the points of concentration $a_i$.
  \begin{pro}\label{p:38w}  Let $\e$ be small enough. For each $i \leq q$,  it holds:
  $$  - \frac{c_2}{2} \sum_{j\leq q; j\neq i} \a _j  \frac{1}{\l _i}\frac{\partial\e _{ij}}{\partial a_i } -  \frac{ \a _i} {K(a_i)}\frac{c_6}{\l _i}{\n K_1}(a_i) = O \Big( \frac{1}{\l_i^2} + \sum_{j > q} \e_{ij} + R_{a_i}  + \sum \frac{\ln \l_k }{\l_k ^{n/2}} \Big) . $$
  Furthermore, for each $ i \geq q+1$, it holds
$$  \frac{|\nabla K(a_i)|}{\l_i} \leq c \Big(\frac{1}{\l_i^3}+ \frac{1}{(\l_id_i)^{n-2}} + \sum_{j\neq i}\e_{ij} + R^2(\e,a,\l) + \sum \frac{\ln \l_k }{\l_k ^{n/2}} \Big) $$
  where   $ R_{a_i} $ and $ c_6 $ are defined in Proposition \ref{p:38} and $R(\e,a,\l)$ is defined in  Proposition \ref{eovv}.
\end{pro}





\subsection{ Proofs of Theorems \ref{th:t2}, \ref{th:t3}} We start by providing the proof of Theorem \ref{th:t2}.

\begin{pfn}{ \bf of Theorem \ref{th:t2}}
Recall that the proof of Theorem \ref{th:t1} relies on the Eqs $(E_i)$, $(F_i)$ and \eqref{alpha1} which follow from Propositions  \ref{p:35}, \ref{p:33} and \ref{p:38}. When $\o \neq 0$, the counterpart of these propositions are Propositions \ref{p:35w}, \ref{p:33w} and \ref{p:38w}. In the new propositions, the new terms are $O(1/\l_i ^{(n-2)/2})$ which is  $O(1/\l_i ^{5/2} )$ if $n \geq 7$ and therefore it is small with respect to the other principal terms. Hence, the Eqs  $(E_i)$, $(F_i)$ and \eqref{alpha1} are not changed. Thus the proof can be repeated exactly by the same way.
\end{pfn}

 \begin{pfn}{ \bf of Theorem \ref{th:t3}}
 We notice that, in the proof of Theorem \ref{th:t1}, for $a_i\in \mathbb{S}^n_+$, the term $\frac{\D K(a_i)}{\l_i^2}$ is a principal term. However, for the dimension $5$, when $\o \neq 0$, another term appears in Proposition \ref{p:33w} which is $\frac{\o(a_i)}{\l_i^{3/2}}$. This term will dominate the previous one. Hence, we will have a changement in the behavior analysis.

We order all the $\l_i$'s, for $i=1,\cdots,N$: $\l_1 \leq \cdots \leq \l_N$ and as in the previous case, we denote by
 $$ I_b := \{ i : a_i \in \partial \mathbb{S}^n_+ \} \qquad \mbox{ and } \qquad I_{in} := \{ i : a_i \in  \mathbb{S}^n_+ \} ,$$
 $$ I' := \{ i: \l_i /\l_1 \to \infty \} \qquad \mbox{ and } \qquad I:= \{1, \cdots , N\} \setminus I'.$$
Multiplying the equation of Proposition \ref{p:33w} by $2^i \times M$ if $ i>q$ (with $M$ a large constant to dominate the $O(\e_{ij})$ which appears in Proposition \ref{p:33w} for $i \leq q$) and $ 2^i$ if $ i\leq q$, and summing over $i$, we obtain:
\be\label{*1*}
\e + \sum_{\i \neq j} \e_{ij} + \sum_{i > q}\frac{1}{(\l_i d_i)^3} = O(\frac{1}{\l_1}).\ee
  This implies that $R_{1,\l}$ and $R_{2,\l}$ defined in Proposition \ref{p:33w} and $R_{a_i}$ defined in Proposition \ref{p:38w} satisfy
  \be\label{*2*} | R_{1,\l} | + | R_{2,\l} |   \leq c \ln (\l_1) /{\l_1^{5/3}} \qquad \mbox{ and } \qquad | R_{a_i} | \leq c / {\l_1^{4/3}} \quad \mbox{ for each } i .\ee

 \begin{lem}\label{wwww1}
$(1)$  Assume that $I_{in} \neq \emptyset$, then $$ \e + \sum_{ i \in I_{in} ; k\neq i }\e_{ik} + \sum_{i \in  I_{in}}\frac{1}{(\l_i d_i)^3}  + \sum_{i \in I_{in}} \frac{c}{\l_i^{3/2}} = O\Big( R_{2,\l} + \sum \frac{\ln\l_k}{\l_k^{5/2}}\Big) = O\Big( \frac{\ln\l_1 }{\l_1^{5/3}} \Big).$$
This implies that $I \cap I_{in} = \emptyset$.\\
$(2)$ For each $k \in I_b$, it holds: $$ \e + \sum_{ i \in I_{b} ; i\neq k }\e_{ik} = O\Big(  \frac{1}{\l_k } + \frac{\ln\l_1}{\l_1^{5/3}}\Big).$$
$(3)$ For each $i \in I_b $ and $k \in I_b$, it holds: $$ ({1}/{\l_i}) |{\partial \e_{ik}}/{\partial a_i} | = o( {1}/{\l_1}).$$
 \end{lem}
 \begin{pf} Multiplying the second equality in  Proposition \ref{p:33w} by $2^i$ and summing over $i \in I_{in}$ the first claim follows.  Using Claim $(1)$  the second one follows by multiplying the first equality in Proposition \ref{p:33w} by $2^i$ and summing over $i \geq k$ and $i \in I_b$.\\
  Concerning the last one, if $i , k \in I$ then, we derive using \eqref{*1*} that $${(1}/{\l_i} ) |{\partial \e_{ik}}/ {\partial a_i} |  \leq c \e_{ik}^{(n-1)/(n-2)} \leq c/ \l_1 ^{(n-1)/(n-2)}.$$
 In the other case, using Claim (2), we obtain ${(1}/{\l_i} ) |{\partial \e_{ik}}/ {\partial a_i} |  \leq c \e_{ik} = o(1/\l_1)$.
 \end{pf}

 Next we claim that :\\
 {\bf Claim A:} $ I = I_b$ and $I_{in} = \emptyset$.\\
 To prove this claim, we follow the proof of Proposition \ref{comparable}. In fact, arguing by contradiction, assume that either $I \neq I_b$ or $I_{in} \neq \emptyset$. Then, from Lemma \ref{wwww1}, it follows that $\e = o(1/\l_1)$ (which gives Claim 1). Claim 2 follows  immediately from (1) of Lemma \ref{wwww1}. Claim 3 follows from Proposition \ref{p:38w} and $(3)$ of Lemma \ref{wwww1}. Claim 4 and $(i)$ of Claim 5 follow from Proposition \ref{p:33w} and Lemma \ref{wwww1}. $(ii)$ of Claim 5 follows from $(i)$ and the fact that all the $\e_{ij}$'s, for $i,j \in A_z$, are of the same order. Finally, the sequel of the proof follows exactly in the same way which completes the proof of Claim A.

\medskip
 Now, we are in the same situation with Theorem \ref{th:t1} but with $\ell = 0$ and the sequel of the proof can be repeated exactly in the same way.
\end{pfn}

\section{ Proof of Theorem \ref{th:t4}}

Arguing by contradiction we assume that there exists  a sequence $(u_k)$  of energy bounded solutions of $(\mathcal{P}_0)$ which blows up. Thus two cases may occur:  $(1)$  either $(u_k) $ converges weakly to $0$, $(2)$  or $(u_k) $ converges weakly to a positive solution $\o$ of $(\mathcal{P}_0)$.
In both  cases, we will run to  a contradiction.\\
We start by considering  the first case, that is when  $(u_k) $ converges weakly to $0$.\\
In this case, $(u_k)$ has to blow up and therefore it will enter in some $V(q,\ell,\tau)$. Arguing as in the proof of Theorem \ref{th:t1} and using the same notation. Observe that Propositions \ref{eovv}, \ref{p:35}, \ref{p:33}, \ref{p:38} hold true with $\e = 0$. Furthermore, Lemmas \ref{sum}, \ref{56}, \ref{55} and \ref{57} hold  with $\e = 0$. Hence, in this case,   for each $i = 1,\cdots, q + \ell$, Eqs $(E_i')$ and $(F_i')$ (defined in \eqref{Ei'} and \eqref{Fi'}) become

\be \label{tildeEi'} (\tilde{E}_i') : \, \,  \begin{cases}\displaystyle{  \frac{ \D K(a_i) }{\l_i ^2 K(a_i)} = o\big( \frac{1}{\mu_1} \big) } & \displaystyle{\mbox{for } i\in I_{in}, } \\
\displaystyle{    - \frac{c_2} {2}\sum _{ i \neq j\in I_b} \a _j \l _i\frac{\partial \e _{ij}}{\partial \l _i}  -  \frac{\a_i }{K(a_i)}  \frac{c_3}{\l _i}\frac{\partial K}{\partial \nu}(a_i)  = o\big( \frac{1}{\mu_1^{(n-1)/(n-2)}} \big) } & \displaystyle{\mbox{for } i\in I_{b} } ,\end{cases}  \ee
 \be\label{tildeFi'}  (\tilde{F}_i') : \quad  - \frac{c_2}{2} \sum_{i \neq j\in I_b} \a _j  \frac{1}{\l _i}\frac{\partial\e _{ij}}{\partial a_i } -  \frac{ \a _i} {K(a_i)}\frac{c_6}{\l _i}{\n K_1}(a_i) = o \big( \frac{1}{\mu_1^{(n-1)/(n-2)} } \big)  \, \, \,  \mbox{for } i\in I_{b} . \ee

Now we claim that:

{\bf Claim 1:} For each $i\in I_{in}$, it holds that $\mu _i / \mu_1 \to \infty$ as $k \to + \infty$, (that is $I\cap I_{in} = \emptyset$).\\
In fact, arguing by contradiction, assume that there exists $i \in I \cap I_{in}$. Thus, using Lemma \ref{56}, we get that $a_i$ has to converge to a critical point $y\in \mathbb{S}^n_+$ of $K$ and therefore $|\D K(a_i)| \geq c > 0$ for $k $ large. Thus Eq $(\tilde{E}_i')$ gives a contradiction and therefore the claim follows.

\medskip

Now summing $2^i (\tilde{E}_i')$ for $i \in I_b \setminus I$ we derive that
\be \label{emna17} \sum_{ i \in I_b \setminus I \, ; \, i \neq k \in I_b  } \e_{ik} = o\Big( \frac{1}{\mu_1}\Big). \ee

From Eq $(\tilde{F}_1')$ and Lemma \ref{57} we derive that $a_1$ has to converge to a critical point $z_1$ in $\partial \mathbb{S}^n_+$. Using Eq $(\tilde{E}_1')$ and \eqref{emna17}   we get that
\begin{equation} \label{kkk0}\sum_ {j \neq 1; j\in I_b} \e_{1j} = \frac{ c \, \,  (1+o(1)) }{  \l_1 } . \end{equation}
 Hence, let $j$ be such that $ d(a_1,a_j) := \min \{ d(a_1,a_k): k \in I \cap I_b \}$, then it holds that $d( a_j , a_1 ) \to 0$ which implies that $\# B_{z_1} \geq 2$ where $B_{z_1} := \{ i\in I \cap I_b : a_i \to z_1 \}$.

 We introduce the following sets: $$ A'_{z_1}:= \{ k \in B_{z_1}: \lim d( a_k , a_{1} ) / d( a_{j } , a_{1} ) = \infty\} \quad \mbox{ and } \quad A_{z_1} := B_{z_1} \setminus A'_{z_1}.$$

{\bf Claim 2:}  For each $i\neq k \in A_{z_1}$, it holds $\l_i ^{(n-3)/(n-2)} d(a_i,a_k) $ is bounded above and below.

We notice that, for each $i\neq k \in A_{z_1}$, it holds that $d(a_i,a_k)$ and $d(a_{j}, a_{1})$ are of the same order. Furthermore, $\l_i$ and $\l_k$ are of the same order. Hence, $\e_{ik} = (\l_i \l_k d(a_i,a_k)^2)^{(2-n)/2} (c +o(1))$ and all the $\e_{ik}$, for each $i\neq k \in A_{z_1}$, are of the same order. Thus Claim 2  follows immediately from \eqref{kkk0}.

\medskip

To conclude the proof of the theorem, as the proof of Proposition \ref{comparable}, we need to multiply $(\tilde{F}_i')$ by $\a_i \l_i (\ov{a} -\langle a_i,\ov{a}\rangle a_i)$ and summing for $i\in A_{z_1}$ where $\ov{a}$ is the barycenter of the points $a_i$'s for $i\in A_{z_1}$. Observe that Claim 2 implies that $\l_i | \ov{a} -\langle a_i,\ov{a}\rangle a_i | \leq c \l_i d(a_i, \ov{a} ) \leq c \l_i ^{1/(n-2)}$.
It holds that
\begin{equation} \label{kkk1} c_6 \sum_{i\in A_z} \a_i^2 \frac{\n K_1(a_i) }{K_1(a_i)} (\ov{a} -\langle a_i,\ov{a}\rangle a_i)  - \frac{c_2}{2} \sum_{i, j\in A_z; j\neq i} \a_i \a_j \frac{\partial \e_{ij}}{\partial a_i}  (\ov{a} -\langle a_i,\ov{a}\rangle a_i) = o\big( \frac{1}{\mu_1} \big).\end{equation}
Observe that,  \eqref{alphai0}  implies $ \a_i^2 = K(a_i)^{(2-n)/2} + o(1)$. Furthermore, using Lemmas \ref{lemderepsilon} and \ref{derK}, there hold
\begin{align}
& \sum_{i\in A_{z_1}} \frac{ 1 }{K_1(a_i)^\frac{n}{2} }  \n K_1(a_i)  (\ov{a} -\langle a_i,\ov{a}\rangle a_i)  = \sum_{i\in A_{z_1}} O( | \ov{a} - a_i |^2) = O\Big(\frac{1}{\mu_1^\frac{2(n-3)}{(n-2)}}\Big) = o\Big(\frac{1}{\mu_1}\Big), \label{kkk2} \\
& \frac{\partial \e_{ij}}{\partial a_i}  (\ov{a} -\langle a_i,\ov{a}\rangle a_i) + \frac{\partial \e_{ij}}{\partial a_j}  (\ov{a} -\langle a_j,\ov{a}\rangle a_j) \geq  c \e_{ij}\label{kkk3} . \end{align}
Eqs \eqref{kkk2} and \eqref{kkk3} imply that Eq \eqref{kkk1} cannot occur. Hence the first case (that is $(u_k) $ converges weakly to $0$) cannot occur.\\
Now we will focus on  the second case: that is $(u_k) $ converges weakly to $\o \neq 0$.\\
The proof follows the previous one but we need to use Section 4 instead of Section 3. Observe that when $\o \neq 0$, comparing with Propositions \ref{eovv}-\ref{p:38}, the new propositions contain other terms due to the presence of $\o$.\\
For $n \geq 7$, these terms are small with respect to the other principal terms (which exist in Propositions \ref{eovv}-\ref{p:38}). Hence the new propositions seem as there exist no changes and therefore the previous proof work exactly in the same way.\\
However, for $n=5$, the new term $\o(a_i)/\l_i^{3/2}$ dominates $\D K(a_i)/\l_i ^2$ and therefore the principal terms in Proposition \ref{p:33w} will be different comparing with Proposition \ref{p:33}. \\
Recall that Section 4 holds for $\e =0$. Using the same notation than the proof of Theorem \ref{th:t3}, we derive that \eqref{*1*} and \eqref{*2*} hold with $\e =0$.  Now, from the first assertion of Lemma \ref{wwww1}, we derive that $a_1 \in \partial\mathbb{S}^5_+$. Moreover, using Proposition \ref{p:38w} and Lemma \ref{wwww1}, we obtain that $a_1$ converges to a critical point $z$ of $K_1$. In addition, using Proposition \ref{p:33w} we get that
\be \label{z*1} \sum \e_{1k} = c \, \frac{\partial_\nu K(z)}{\l_1} (1+o(1))\ee
which implies that  $\partial_\nu K(z)$ has to be positive (if not the previous equality cannot occur).
Observe that \eqref{z*1} implies that there exists at least one index $j\in I$ such that $d(a_j,a_1) \to 0$ as $k \to \infty$.\\
As before, let $B_z:= \{ j \in I: a_j \to z\}$. It holds that $\# B_z \geq 2$. Furthermore, following the end of the proof of Proposition \ref{comparable}, we deduce that \eqref{z*2} holds. Finally observe that \eqref{z*3} and \eqref{z*4} imply that \eqref{z*2} cannot occur. Thus the proof of our theorem is complete.


\section{Appendix}

 In this section we collect some technical Lemmas used in this paper.

\begin{lem}\cite{AB20a} \label{lem:varphi}
For $a\in \partial \mathbb{S}_+^n$, we have $\partial \delta_{a,\l} /\partial \nu = 0$ and therefore $\varphi_{a,\l} = \d_{a,\l}$. For $a\notin \partial \mathbb{S}_+^n$, we have
$$ \d_{a,\l} \leq  \varphi_{a,\l} \leq 2 \d_{a,\l}\, \, ; \quad | \l \partial \varphi_{a,\l} /\partial \l | \leq c \d_{a,\l} \, \, ; \quad  | (1/\l) \partial \varphi_{a,\l} /\partial a^k | \leq c \d_{a,\l}, \leqno{(i)}$$
where $a^k$ denotes the $k$-th component of $a$.
$$   \varphi_{a,\l}  = \d_{a,\l} + c_0 \frac{H(a,.)}{\l^{(n-2)/2}} + f_{a,\l}\qquad \mbox{ where } \leqno{(ii)} $$
$$ | f_{a,\l} | _{\infty} \leq  \frac{c}{(\l d_a)^2} \frac{H(a,.)}{\l^\frac{n-2}{2}} \leq  \frac{c}{\l^\frac{n+2}{2} d_a^{n}} \, \, ; \, \,  | \l \frac{\partial f_{a,\l} }{\partial \l}  | _{\infty} \leq \frac{c}{\l^\frac{n+2}{2} d_a^{n}} \mbox{ and }  | \frac{1}{\l} \frac{\partial f_{a,\l} }{\partial a^k}  | _{\infty} \leq \frac{c}{\l^\frac{n+4}{2} d_a^{n+1}},$$ where $d_a:= d(a,\partial \mathbb{S}^n_+$). Furthermore, it holds that ${H(a,.)}/{\l^{(n-2)/2}} \leq c \d_{a,\l}$.
 \end{lem}

 \begin{lem}\label{lowerL2}\cite{AB20b}
 Let $a\in \ov{ \mathbb{S}^n_+}$ and $\l > 0$ be large.\\
  $(i)$ Assume that $\e \ln \l$ is small enough, then it holds
 \begin{align} \label{lower2}\d_{a,\l}^{-\e}(x) = & c_0^{-\e} \l^{-\e (n-2)/2}\Big(1 + \frac{n-2}{2}\, \e \ln (2+(\l^2-1)(1 -\cos d(a,x))) \Big) \\
 & + O\Big(\e ^2 \ln (2+(\l^2-1)(1 -\cos d(a,x)))\Big) \quad \mbox{ for each } y \in \ov{\mathbb{S}^n_+}. \nonumber \end{align}
 $(ii)$ For each $\g > 0$ and each $\b\in [0, n/(n-2))$, it holds
 $$ 0 < \int_{\mathbb{S}^n_+} \d_{a,\l}^{p+1 - \b } (x) \ln^\g \Big(2+(\l^2-1)(1 -\cos d(a,x))\Big)dx = O\Big( \frac{1}{\l^{\b (n-2)/2} }\Big).$$
 \end{lem}

\begin{lem}\label{lemderepsilon} \cite{AB20a}
Let $i,j \in I_b$ be such that $\l_i$ and $\l_j$ are of the same order and $| a_k - h | \to 0$ for $k=i,j$ for some $h \in \partial \mathbb{S}^n_+$. Then we have
$$ e_{ij}:=  \frac{ \partial \e_{ij} }{\partial a_i } (h- \langle a_i , h\rangle a_i) + \frac{ \partial \e_{ij} }{\partial a_j } (h- \langle a_j , h\rangle a_j)  \geq  c \, \e_{ij}. $$
\end{lem}

\begin{lem}\label{derK}
$(1)$ Let $a, \, h \in \partial \mathbb{S}^n_+$ be close to a critical point  $z$ of $K_1$. Then it  holds that
\be\label{mm4} \frac{1}{K_1(a)^{n/2}} \n K_1(a) \big(h - \langle a, h \rangle a\big) =  - \frac{1}{K_1(h)^{n/2}} \n K_1(h) \big(a - \langle a, h \rangle h\big)  + O(  | a - h |^2 ). \ee
$(2)$ Let $e\in \partial \mathbb{S}^n_+$ be such that $|e - z| \geq c > 0$.
$$ \frac{1}{K_1(a)} \n K_1(a) \big(e - \langle a, e \rangle a\big) =  \frac{1}{  K_1(z)} D^2 K_1(z) \big(a - \langle a, z \rangle z, e - \langle e, z \rangle z\big) + O(  | a- z |^2 ). $$
\end{lem}
\begin{pf}
Let $$ \beta(t):= \frac{ h + t(a-h)}{| h + t(a-h) |} \quad , \quad g(t):= \frac{ 2/(n-2) }{ K_1(\beta(t))^{(n-2)/2}} \quad \mbox{ for } t\in [0,1].$$
It is easy to get that
$$ \beta'(t) = \frac{1} { | h + t(a-h) |} \Big( a - h - \langle \beta(t) , a-h \rangle \beta(t) \Big) \quad , \quad  \langle \beta(t) , a-h \rangle = O(| a - h |^2),$$
(since we have $| a- h |^2 = 2(1-\langle a,h\rangle)$) and therefore it holds that $| \beta'(t) | = |a-h | (1+o(1))$ uniformly in $t\in [0,1]$. Furthermore, easy computations imply that $| \beta ''(t) | = O( |a-h | ^2)$ uniformly in $t\in [0,1]$. On the other hand, we have
$$ g'(t) = \frac{-1}{  K_1(\beta(t))^{n/2}} \n K_1(\beta(t)) \big(\beta'(t)\big) $$
 and, since $a$ and $h$ are close to a critical point $z$ of $K_1$, we derive that
\be\label{mm3}g''(t) = o( | \beta'(t) |^2 ) - \frac{1}{  K_1(\beta(t))^{n/2}} D^2 K_1(\beta(t)) \big(\beta'(t), \beta'(t)\big) + o( | \beta''(t) |  )= O ( | a- h |^2 ) \ee(uniformly in  $ t\in [0,1]$). Now,
\begin{align*} \frac{1}{K_1(a)^{n/2}} \n K_1(a) \big(h - \langle a, h \rangle a\big)  & + \frac{1}{K_1(h)^{n/2}} \n K_1(h) \big(a - \langle a, h \rangle h\big) \\
 & = g'(1) - g'(0) = \int _0^1 g''(t) \, dt \end{align*}
  which implies the first claim of the  lemma.\\
To prove the second claim,  since $z$ is a critical point of $K_1$, it follows that $1/K_1(a) = 1/K_1(z) + O(| a - z |^2)$. Furthermore, let $f(a):= \n K_1(a) \big(e - \langle e, a \rangle a\big)$. It follows that: $f(a)=f(z) + f'(z) ( a - \langle a, z \rangle z) +  O(| a - z |^2)$ which completes the proof.
\end{pf}

\begin{rem} Let $m=3$ and $z$ be a critical point of $K_1$. Assume that the matrix associated to $D^2K_1(z)$ has a positive eigenvalue $\s >0$. Then the function $\mathcal{F}_{z,3}$ has a critical  point  of the form $(\ov{b},0,-\ov{b})$.

  In fact, each critical point of $\mathcal{F}_{z,3}$ has to satisfy the following system

$$\begin{cases} & D^2 K_1(z) (x_1,.) - (n-2) \frac{x_1 - x_2}{| x_1 - x_2|^n} - (n-2) \frac{x_1 - x_3}{| x_1 - x_3|^n} = 0\\
& D^2 K_1(z) (x_2,.) - (n-2) \frac{x_2 - x_1}{| x_2 - x_1|^n} - (n-2) \frac{x_2 - x_3}{| x_2 - x_3|^n} = 0\\
& D^2 K_1(z) (x_3,.) - (n-2) \frac{x_3 - x_1}{| x_3 - x_1|^n} - (n-2) \frac{x_3 - x_2}{| x_3 - x_2|^n} = 0
\end{cases} $$
Taking $(x_1,x_2,x_3) = (\ov{b},0,-\ov{b})$ then the second equation is  satisfied. The first and the third equations give the same one which is : $$ D^2 K_1(z) (\ov{b},.) - (n-2) \frac{\ov{b}}{| \ov{b} |^n} - (n-2) \frac{2 \ov{b}}{|2 \ov{b} |^n} =  D^2 K_1(z) (\ov{b},.) - (n-2) \frac{\ov{b}}{| \ov{b} |^n} (1 + \frac{1}{2^{n-1}}) = 0 .$$
This equation has a solution $\ov{b}= \g e_\s$ where $e_\s$ is an eigenvector (with norm 1) associated to the eigenvalue  $\s > 0$ and $\g> 0$ satisfies $$ \g^n = (n-2) (1 + \frac{1}{2^{n-1}})\frac{1}{\s}.$$

This remark implies that we cannot prove that $\l_i ^{(n-2)/n} | a_i - z | \geq c $ for each $i$ (see  Lemma \ref{nonsimple1}) since it is possible that for one index $j$ we have  $\l_j ^{(n-2)/n} | a_j - z | \to 0 $. But this can occur for at most one index.
 \end{rem}

\begin{lem}\label{somest} For each $k\neq j$, It holds
$$ \int (\d_k \d_j)^\frac{n}{n-2} \leq c \e_{kj} ^\frac{n}{n-2} \ln\e_{kj}^{-1}  \quad \& \quad   \int \d_k \d_j^\frac{4}{n-2} \leq c \e_{kj} ^\frac{n}{n-2} \ln\e_{kj}^{-1} + c \frac{\ln \l_k}{\l_k^{n/2}} + c  \frac{\ln \l_j}{\l_j^{n/2}} .$$
\end{lem}

\begin{pf} The first assertion is extracted from \cite{B1} (Estimate 2 page 4). For the second one,  let $B_i := B(a_i,1)$, we have
\begin{align*}   &  \int \d_k \d_j^\frac{4}{n-2} \\ & \leq \int _{[\d_j \leq \d_k]  \cap B_k } \d_k (\d_j\d_k)^\frac{2}{n-2} +  \int _{[\d_j \leq \d_k]  \cap B_k ^c}  \d_k^\frac{n+2}{n-2}  +  \int _{ [\d_k\leq \d_j] \cap B_j }( \d_k \d_j)^\frac{2}{n-2} \d_j +  \int _{[\d_k\leq \d_j]  \cap B_j ^c} \d_j^\frac{n+2}{n-2} \\
& \leq \Big( \int (\d_j\d_k)^\frac{n}{n-2}\Big)^\frac{2}{n}  \Big( \Big(\int _{B_k}\d_k ^\frac{n}{n-2} \Big)^\frac{n-2}{n} + \Big(\int _{B_j}\d_j ^\frac{n}{n-2} \Big)^\frac{n-2}{n} \Big) + \frac{ c }{ \l_k ^{(n+2)/2}} + \frac{ c }{ \l_j ^{(n+2)/2}} \\
& \leq c\Big( \e_{kj} ^\frac{n}{n-2} \ln\e_{kj}^{-1}\Big)^\frac{2}{n} \Big( \Big(\frac{ \ln \l_k}{ \l_k ^{n/2}}\Big)^\frac{n-2}{n} + \Big(\frac{ \ln \l_j}{ \l_j ^{n/2}}\Big)^\frac{n-2}{n} \Big) + \frac{ c }{ \l_k ^{(n+2)/2}} + \frac{ c }{ \l_j ^{(n+2)/2}} \\
& \leq c  \e_{kj} ^\frac{n}{n-2} \ln\e_{kj}^{-1} + c \frac{ \ln \l_k}{ \l_k ^{n/2}} + c \frac{ \ln \l_j}{ \l_j ^{n/2}}.
\end{align*}
Hence the proof  is completed.  \end{pf}

\begin{lem}\label{Qwalpositive} Assume that the $\e_{ij}$'s are small enough. Let $v \in E_{\o,a,\l} ^\perp$ (defined in \eqref{Ewal}). Written $v$ as
$$ v := v_- + v_0 + v_+ \quad  \mbox{ with } \quad
 v_- \in N_-(\o) \, \, ; \, \,  v_0 \in H_0(\o) \, \, ; \, \, v_+ \in N_+(\o) $$ { where }  $ H_0(\o):= span(\o) \oplus N_0(\o)$,  (defined in \eqref{vvv*1}). Then there exists a positive constant $\underline{c} $ such that
\begin{align}  & \| v_0 \| = o( \| v \|  ),  \label{0q0} \\
& Q_{\o,a,\l} (v_- ) = Q_\o (v_-) + o( \| v_- \| ^2 ) \leq -\,  \underline{c}\,  \| v_-\|^2, \label{0q-}\\
& Q_{\o,a,\l} (v_+ ) \geq  \underline{c} \,  \| v_+\|^2 + o( \| v \| ^2), \label{0q+}
\end{align}
where $Q_\o$ is defined in  \eqref{qw} and $Q_{\o,a,\l}$ is defined by
\be \label{ff7} Q_{\o,a,\l} (v) = \| v \|^2 - p \sum \int_{\mathbb{S}^n_+} \d_i ^{p-1} v ^2 - p  \int_{\mathbb{S}^n_+} K \o ^{p-1} v ^2. \ee

This implies that the quadratic form $Q_{\o,a,\l}$ is a non degenerate one in the space $E_{\o,a,\l} ^\perp$.
\end{lem}

\begin{pf} First note that  the spaces  $H_0(\o)$, $N_-(\o)$ and $N_+(\o)$  are orthogonal spaces with respect to $\langle.,.\rangle$ and the associated bilinear form $B_\o(.,.)$ ($:= \int_{\mathbb{S}^n_+}  K \o^{p-1} . . $).

We start by proving \eqref{0q0}. Since $v_0 \in H_0(\o)$, it follows that $ v_0 = \g_0 \o + \sum \g_i e_i$ where $(e_1,\cdots,e_m)$ is an orthonormal basis of $N_0(\o)$. Using the fact that $v \in E_{\o,a,\l}^\perp$ (which implies that $ v \perp u_{\a,\b}$ and $v \perp \partial u_{\a,\b} / \partial \b_i$ for each $i $), it follows that
\begin{align*} & \g_0 =  \langle v_0 , \o \rangle  =  \langle v , \o \rangle = (1/ \a) \langle v , u_{\a,\b}  \rangle  -  (1/\a) \sum \b_i \langle v , e_i \rangle = o( \| v \| ), \\
& \g_i = \langle v_0 , e_i \rangle = \langle v , e_i \rangle = \langle v , e_i + (\partial h(\b)/ \partial \b_i)  \rangle - \langle v ,  \partial h(\b)/ \partial \b_i  \rangle = o( \| v \| ) \quad \forall \,\,   1 \leq i \leq m \end{align*}
(by using the smallness of $\b$ and $h(\b)$ in the $C^1$ sense with respect to $\b$). This ends the proof of \eqref{0q0}.

Concerning \eqref{0q-}, we have
$$ Q_{\o,a,\l} (v_- ) = \| v_- \| ^2 - p \sum \int \d_i ^{p-1} v_- ^2 - p \int K \o ^{p-1} v_- ^2 = Q_\o ( v_-) - p \sum \int \d_i ^{p-1} v_- ^2 .$$
Observe that, since $v_-$ belongs to a  fixed finite dimensional space, we derive that $\| v_- \|_{\infty} \leq c \| v_- \| $ and therefore
$$ \int \d_i ^{p-1} v_-^2 \leq \| v_- \|^2_{\infty}  \int \d_i ^{p-1}  = o(  \| v_- \|^2 ) \quad \mbox{ for each } i.$$
Hence the proof of \eqref{0q-} follows by using \eqref{qw+-}. \\
It remains to prove \eqref{0q+}.
Note that, using Proposition 3.1 of \cite{B1}, there exists a constant $\underline{c}_1> 0$ such that
\begin{equation}\label{h1}
\| h \| ^2 - \frac{n+2}{n-2}  \sum_{i=1}^N \int _{\mathbb{S}^n_+} \d_{a_i, \l_i} ^{\frac{4}{n-2}} h^2  \geq \underline{c}_1 \| h \| ^2 \quad \mbox{ for each } h \in E_{a,\l}^\perp.
\end{equation} where $E_{a,\l}^\perp$ is introduced in Proposition \ref{lambdaepsilon}. \\
In addition, the sequence of the eigenvalues (denoted by $(\s_i)$)  corresponding to $Q_\o$ (defined by \eqref{qw}) satisfies $\s_i \nearrow 1$. Let $N_k(\o)$ be the eigenspace associated to the eigenvalue $\s_k$. These spaces are orthogonal with respect to $\langle ., .\rangle$ and the bilinear form $B_\o$. Let $\s_{k_0}:= \min\{ \s_i : \s_i > 0\}$. Hence it is easy to see that $N_+(\o)= \oplus_{ k \geq k_0} N_k(\o)$. Furthermore, it holds
\begin{equation}\label{h2} Q_\o (h) \geq \s_k \| h \| ^2 \qquad \mbox{ for each } h \in \oplus_{ j \geq k } N_j(\o).\end{equation}
Let $k_1$ be such that $\s_{k_1} \geq 1- \underline{c}_1 /4 $. We decompose $N_+(\o)$ as follows:
$$ N_+(\o) := (\oplus_{k_0 \leq k \leq k_1} N_k(\o)) \oplus (\oplus_{k > k_1} N_k(\o):= N_+^{0,1}(\o) \oplus N_+^{1}(\o).$$
Note that $ N_+^{0,1}(\o)$ is a fixed finite dimensional  space.
Now, since  $v_+ \in  N_+(\o)$, it holds that
\be\label{qaz2} v_+:= v_0^+ + v_1^+ \qquad \mbox{ where } v_0^+ \in  N_+^{0,1}(\o) \mbox{ and } v_1^+ \in  N_+^1(\o).\ee
Hence, using the orthogonality of the spaces $N_+^{0,1}(\o)$ and $ N_+^1(\o)$, it follows that
\begin{align*} Q_{\o,a,\l} (v_+) = \| v_0^+ \| ^2 +  \| v_1^+ \| ^2 & - p  \sum_{i=1}^N \int _{\mathbb{S}^n_+} \d_{a_i, \l_i} ^{\frac{4}{n-2}} \{ (v_0^+)^2 + (v_1^+)^2 + 2 v_0^+ v_1^+ \} \\
 & - p  \int _{\mathbb{S}^n_+} K \o ^{\frac{4}{n-2}} \{ (v_0^+) ^2 + (v_1^+)^2 \}. \end{align*}
Observe that
\begin{align*}
& \| v_0^+ \| ^2 - \frac{n+2}{n-2}  \int _{\mathbb{S}^n_+} K \o ^{\frac{4}{n-2}} (v_0^+) ^2 = Q_\o(v_0^+)  \geq \s_{k_0} \| v _0^+ \| ^2  \quad \mbox{ (by using \eqref{h2}) } \\
&  \int _{\mathbb{S}^n_+} \d_{a_i, \l_i} ^{\frac{4}{n-2}} (v_0^+)^2 \leq \| v_0^+\|^2_\infty \int _{\mathbb{S}^n_+} \d_{a_i, \l_i} ^{\frac{4}{n-2}}  = o( \| v_0^+ \| ^2)\\
&  \int _{\mathbb{S}^n_+} \d_{a_i, \l_i} ^{\frac{4}{n-2}} | v_0^+ | | v_1^+ | \leq \| v_0^+\|_\infty \int _{\mathbb{S}^n_+} \d_{a_i, \l_i} ^{\frac{4}{n-2}}  | v_1^+ | =  o( \| v_0^+ \| \| v_1^+ \|) = o( \| v_1^+  \| ^2 + \| v_0^+ \|^2 )\\
& - p  \int _{\mathbb{S}^n_+} K \o ^{\frac{4}{n-2}}  (v_1^+)^2 = Q_\o (v_1^+) - \| v_1^+ \|^2 \geq (\s_{k_1} - 1 ) \| v_1^+ \| ^2 \geq -(\underline{c}_1/4) \| v_1^+ \| ^2
\end{align*}
where, for the last formula,  we have used  \eqref{h2} and the choose of $k_1$. Combining these estimates, we get
\be\label{qaz1} Q _{\o,a,\l} (v_+)  \geq \{ \| v_1^+ \| ^2  - p  \sum_{i=1}^N \int _{\mathbb{S}^n_+} \d_{a_i, \l_i} ^{\frac{4}{n-2}} (v_1^+)^2 \} + \s_{k_0} \| v _0^+ \| ^2  -\frac{\underline{c}_1}{4} \| v_1^+ \| ^2 + o( \| v_1^+ \| ^2 + \| v_0^+ \| ^2) .\ee
Note that the function $v_1^+$ is not necessarily in $E_{a,\l}^\perp$. For this raison we write $ v_1^+ := \sum t_i \psi _i + \tilde{v}_1^+$ with  $\tilde{v}_1^+ \in E_{a,\l}^\perp$ where the $\psi_i$'s are the functions $ \varphi_j$'s and their derivatives with respect to $ \l_j$ and $a_j^k$. Let $\psi _i \in \{ \varphi_i,\l_i  \partial \varphi_i / \partial \l_i  , (1/\l_i)  \partial \varphi_i / \partial \a_i ^k  \}$, it follows that
\begin{align*} t_i + o(\sum | t_k |) & = \langle v_1^+ , \psi_i \rangle =  \langle v  , \psi_i \rangle -  \langle v_0 , \psi_i \rangle  - \langle v_- , \psi_i \rangle - \langle v_0^+ , \psi_i \rangle\\
&  = O \Big( \int \d_i ^p ( | v_0| + | v_- | + | v_0^+ | \Big) = o(  \| v_0 \| + \| v_- \| + \| v_0^+ \|) = o(  \| v \|)\end{align*} by using the fact that these functions are in  fixed finite dimensional spaces.  Thus we derive that $ \| v_1^+ \|^2 = \| \tilde{v}_1^+ \|^2 +o( \| v \|^2 )$ and therefore, using \eqref{h1}, we get
\begin{align}  \| v_1^+ \| ^2  - p  \sum_{i=1}^N \int _{\mathbb{S}^n_+} \d_{a_i, \l_i} ^{\frac{4}{n-2}} (v_1^+)^2  & = \| \tilde{v}_1^+ \| ^2  - p \sum_{i=1}^N \int _{\mathbb{S}^n_+} \d_{a_i, \l_i} ^{\frac{4}{n-2}} (\tilde{v}_1^+)^2  + o( \| \tilde{v}_1^+ \| ^2  +\| v \| ^2 )\nonumber  \\
& \geq \frac{1}{2} \underline{c}_1  \| v_1^+ \| ^2 + o( \| v \| ^2 )   . \label{qaz3} \end{align}
Combining \eqref{qaz2}, \eqref{qaz1} and \eqref{qaz3}, we get
$$ Q _{\o,a,\l} (v_+)  \geq \frac{1}{4} \underline{c}_1   \|  v_1^+ \| ^2 +  \s_{k_0} \| v _0^+ \| ^2 + o( \| v_1^+ \| ^2 + \| v_0^+ \| ^2) + o( \| v \| ^2 )  \geq c \| v_+ \| ^2 + o( \| v \| ^2 ) .$$
Thus the result follows.
\end{pf}

In the following lemma, we collect some formulae whose proof follows immediately by some standard calculus computations
\begin{lem}\label{ti56} Let $t_i > 0$ and $a,b \in \R$, there hold
\be \label{ti5}
 | (\sum t_i)^\g - \sum t_i ^\g | \leq c \begin{cases} \sum_{i\neq j} (t_it_j)^{\g/2} \quad & \mbox{ if }  0 < \g \leq 2 \\
 \sum_{i\neq j} t_i ^{\g-1} t_j  \quad & \mbox{ if }  \g > 2\end{cases} . \ee
\be \label{ti6}
 | | a+b|^\g - |a| ^\g - \g | a |^{\g-2} a b  | \leq c \begin{cases}  | b | ^\g + | a |^{\g-2} b^2 \quad & \mbox{ if }  \g > 2 \\
 | b | ^\g  \quad & \mbox{ if }  1 < \g \leq 2\end{cases} . \ee
\end{lem}

\end{document}